\documentclass[a4paper,12pt]{amsart}

\usepackage{amsmath}
\usepackage{amsthm}

\usepackage[english]{babel}

\newcommand{\ignore}[1]{}

\newcommand{\From}[0]{From\ }

\newcommand{\Set}[1]{\left\{\, #1 \,\right\}}
\newcommand{\Span}[1]{\left\langle\, #1 \,\right\rangle}

\renewcommand{\phi}[0]{\varphi}
\renewcommand{\theta}[0]{\vartheta}

\newcommand{\F}{\text{$\mathbf{F}$}}

\newcommand{\OldU}[2]{\underbrace{#1 \dots #1}_{#2}}
\newcommand{\U}[2]{#1^{#2}}

\swapnumbers

\numberwithin{dummy}{subsection}

\newtheorem*{lemma}{Lemma}

\newcommand{\uptoasign}{\approx}

\DeclareMathOperator{\End}{End}
\DeclareMathOperator{\ad}{ad}

\begin{document}

%AEC2% According to the instructions of Annali
\bibliographystyle{amsplain}

\author{A.~Caranti and M.~R.~Vaughan-Lee}

\title{Graded Lie algebras of maximal class IV}

\address{Dipartimento di Matematica\\
  Universit\`a degli Studi di Trento\\
  via Sommarive 14\\
  I-38050 Povo (Trento)\\
  Italy} 

\email{caranti@science.unitn.it} 

\urladdr{http://www-math.science.unitn.it/\~{ }caranti/}

\address{Christ Church\\
  University of Oxford\\
  England}

\email{michael.vaughan-lee@christ-church.oxford.ac.uk}

\urladdr{http://users.ox.ac.uk/\~{ }vlee/}

\date{11 July 1999% Version 1.21
  }

\subjclass{17B70 17B65
  17B05 17B30} 

\keywords{Graded Lie algebras of maximal class}

\begin{abstract}
  We describe the  isomorphism classes of certain infinite-dimensional
  graded Lie  algebras of  maximal class, generated  by an  element of
  weight  one  and  an element  of  weight  two,  over fields  of  odd
  characteristic.
\end{abstract}

\thanks{The first author has been partially supported by MURST
  (Italy). The first author is a member of CNR-GNSAGA (Italy), now
  INdAM-GNSAGA. The authors are grateful to CNR-GNSAGA for 
  supporting a visit of the second author to Trento. The second author
  is grateful to the Department of Mathematics of the University of
  Trento for their kind hospitality.}

\maketitle

\thispagestyle{empty}

\section{Introduction}

Let $M$ be a Lie algebra over the field $\F$. Suppose $M$ is
nilpotent of nilpotency class $c$, so that $c$ is the smallest
number such that $M^{c+1}=0$. If $M$ has finite dimension
$n\ge 2$, it is well-known that $c\le n-1$. When $c=n-1$,
$M$ is said to be a Lie algebra of maximal class.

Consider the Lie powers $M^i$. Then $M$ is of maximal
class when the codimension of $M^i$ is exactly $i$,
for $i\le c+1$. It is natural to extend the definition
to an infinite-dimensional Lie algebra $M$ by saying that
$M$ is of maximal class when the codimension of $M^i$ is
$i$ for all $i$ (see \cite{Sha}).

One can grade $M$ with respect to the filtration of the
$M^i$: let 
\begin{equation*}
  L_{i} =  M^{i} /  M^{i+1} ,
\end{equation*}
and consider
\begin{equation}\label{eq:L}
  L = \bigoplus_{i = 1}^{\infty} L_{i}.
\end{equation}
There  is a natural  way of  defining a  Lie product  on $L$,  and the
graded Lie  algebra $L$ has the following  properties: $\dim (L_1)=2$,
$\dim (L_i)\le  1$ for $i\ge 2$,  and $L$ is generated  by $L_1$. Note
that  here too we  allow all  $L_i$ to  be non-zero,  thereby
including
infinite-dimensional  algebras. A  graded Lie  algebra  $L$ satisfying
these  conditions is  called a  graded  Lie algebra  of maximal  class
in~\cite{CMN, CN,  J}. However, this  definition does not  capture all
possibilities. One of the other possibilities for a graded Lie algebra
$L=\bigoplus_{i=1}^\infty L_i$ to be of maximal class is to have $\dim
(L_i)\le 1$ for  all $i\ge 1$, with $L$ generated by  $L_1$ and $L_2$. 
We call
a graded Lie algebra  of this form an  algebra of type
$2$, whereas we refer to
a graded Lie algebra of maximal class in the
sense of~\cite{CMN, CN, J} as an algebra of type $1$.

In studying algebras of type 2, we will mainly deal with
the infinite dimensional ones (as in~\cite{Sha, CMN, CN, J}).
However, our arguments also provide fairly complete
information about finite dimensional algebras.

If the characteristic of the underlying field $\F$ is zero,
it is well-known that there is only one infinite dimensional
algebra of type 1. This is the algebra 
\begin{equation}\label{eq:a}
  a = \Span{x, y : [y \U{x}{i} y] = 0, \ \text{for all $i \ge 1$}},
\end{equation}
where $x$ and $y$ have weight 1. The ideal generated by $y$ is
an abelian maximal ideal here. However, if $\F$ has prime
characteristic $p$ there are uncountably many algebras of
type 1 \cite{Sha, CMN}; these algebras were classified
in \cite{CN, J}.

Over a field $\F$ of characteristic zero there are three
infinite-dimensional algebras of type 2 \cite{SZ, F},
called $m$, $m_2$ and $W$, and these are defined over the
integers. The first one is a close analogue to $a$. It is
given as 
\begin{equation}\label{eq:m}
  m = \Span{e_{1}, e_{2} : 
  [e_{2} \U{e_{1}}{i} e_{2}] = 0, \ \text{for all $i \ge 1$}},
\end{equation}
where $e_1$ has weight 1 and $e_2$ has weight 2. The ideal
generated by $e_2$ is an abelian maximal ideal here.
The second one is defined as 
\begin{equation}\label{eq:m_2}
\begin{aligned}
  m_{2} 
  = 
  \big< \,
  e_{i}, i \ge 1:\
  &[e_{i} e_{1}] = e_{i+1}, \ \text{for all $i \ge 2$},
  \\&[e_{i} e_{2}] = e_{i+2}, \ \text{for all $i \ge 3$}
  \\&[e_{i} e_{j}] = 0, \ \text{for all $i, j \ge 3$}
  \, \big>,
\end{aligned}
\end{equation}
where $e_i$ has weight $i$. Here
$m_2^2=\Span{ e_{i} : i \ge 3}$ is a maximal abelian ideal.
The third algebra is the positive part of the Witt algebra: 
\begin{equation*}
  W = \Span{  e_{i}, i \ge 1: [e_{i} e_{j}] = (i - j) e_{i+j} },
\end{equation*}
and is not soluble.

When one considers these algebras over a field $\F$ of prime
characteristic $p > 2$, ${m}$ and ${m_{2}}$ give algebras of
type 2, but $W$ does not.

We will show in the next section that there is a natural way
to obtain an algebra of type 2 from an uncovered algebra of
type 1. (See the next section for the relevant definition.)
In particular, $m$ arises from $a$ in this way. We will show
that for prime characteristic $p > 2$ the algebras of type 2
consist of
\begin{itemize}
\item algebras arising in this natural way from algebras of type 1,
\item ${m_{2}}$,
\item one further family of soluble algebras,
\item in the case $p=3$, one additional family of soluble algebras.
\end{itemize}
This yields a classification of algebras of type 2 over
fields of characteristic $p > 2$. We believe the case of
characteristic two to be considerably more complicated.

\section{Preliminaries}
\label{sec:preliminaries}

Let $L$ be an \emph{infinite-dimensional}  Lie algebra over a field
$\F$ that is graded over the positive integers:
\begin{equation}\label{eq:L_again}
  L = \bigoplus_{i = 1}^{\infty} L_{i}.
\end{equation}
If  $\dim(L_{1}) =  2$, $\dim(L_{i})  = 1$  for $i  > 1$,
and  $L$ is generated by $L_{1}$, we say that $L$ is \emph{an
algebra of type $1$.} These are the algebras that are called
algebras of maximal class in \cite{CMN,  CN,  J}.
In these papers these algebras are classified over fields
of prime characteristic $p$.

As mentioned in the Introduction, over a field of
characteristic zero there is only one isomorphism class
of algebras of type 1. This is the algebra
$a$ of~\eqref{eq:a} generated by two elements $x$ and $y$ of
weight 1, subject to the relations $[y  \U{x}{i} y] = 0$,
for all $i \ge 1$.  This algebra is  metabelian, and the
graded maximal ideal containing $y$ is abelian. Here we use
the notation
\begin{equation*}
  [y \U{x}{i} y] =  [y \OldU{x}{i} y].
\end{equation*}

If in the algebra~\eqref{eq:L_again} we have $\dim(L_{i}) = 1$
for all $i \ge  1$, and if $L$ is generated by $L_{1}$ and
$L_{2}$, we say that $L$ is  \emph{an algebra  of type $2$.}
Choose non-zero elements  $e_{1} \in L_{1}$ and $e_{2}
\in L_{2}$. Since $L$ is of maximal class, for each $i \ge 2$
we have $[L_{i} L_{1}] = L_{i+1}$. Therefore we can recursively
define $e_{i+1} = [e_{i} e_{1}]$, for $i \ge 2$, and we have
$L_{i} = \Span{e_{i}}$ for all $i$. We keep this notation
fixed for the rest of the paper, allowing ourselves to 
rescale $e_{2}$ when needed.

In \cite{CMN, CN}, to which we refer the reader for all details,
a theory of \emph{constituents} has been developed for
algebras of type $1$ 
over fields of positive characteristic $p$. 
If $L$ is  such an algebra, define its
\emph{$i$-th two-step centralizer} as
\begin{equation*}
  C_{i} = C_{L_{1}}(L_{i}) =
  \Set{ v \in L_{1} : \text{$[u v] = 0$ for $u \in L_{i}$}},
\end{equation*}
for $i > 1$. Each $C_{i}$  is a one-dimensional subspace
of $L_{1}$. A special role is  played by the 
\emph{first two-step centralizer $C_{2}$.}
In fact, the sequence of the two-step
centralizers consists of patterns, called constituents, of
the following type
\begin{equation*}
  C_{i} \ne C_{2}, 
  \quad
  C_{i+1} = C_{i+2} = \dots = C_{i+l} = C_{2},
  \quad
  C_{i+l+1} \ne C_{2}.
\end{equation*}
Here $l$ is  called the \emph{length of the  constituent}.
(We are following the definition  of \cite{CN}, which differs
from that of  \cite{CMN}.) The first constituent requires a
special treatment: its length  is defined as the smallest $f$ 
such
that $C_{f} \ne C_{2}$, and turns  out to be of the form
$f = 2  q$, where $q  = p^{h}$,
for some $h$.  
It is  proved in  \cite{CMN} that 
if the first constituent has length $2 q$, then 
the constituents of $L$ can have lengths of the form
\begin{equation*}
  2 q, 
  \qquad\text{or}\qquad
  2 q - p^{t}, 
  \quad\text{for $0 \le t \le h$.}
\end{equation*}

An algebra  of type 1 is said  to be \emph{uncovered} if  the union of
the  $C_{i}$ does  not  exhaust all  of  $M_{1}$.  
It is proved
in~\cite{CMN} that over  any field  of
positive characteristic there are  uncontably many uncovered
algebras of type  1.  (On the other hand, if the field is
at  most  countable, there  are  algebras of  type  $1$  that are  not
uncovered.)  If $M=  \bigoplus_{i  = 1}^{\infty}  M_{i}$ is  uncovered,
there is an element $z \in M_{1}$ such that
\begin{equation}\label{eq:z}
  \text{$[M_{i} z] = M_{i+1}$ for all $i \ge 1$.}
\end{equation}
We consider the maximal graded subalgebra
\begin{equation*}
  L = \Span{z} \oplus \bigoplus_{i \ge 2} M_{i}
\end{equation*}
of $M$. Because of~\eqref{eq:z}, $L$ is an algebra of
type 2. In addition, the algebra $L$ inherits some kind
of constituent pattern from $M$,
as we will see in the following. 
From now on we will
assume $p > 2$. 

If we apply this procedure to the  unique algebra $M = a$
of~\eqref{eq:a} of type 1 in
characteristic zero,  which is clearly  uncovered, we get  the algebra
$L$ of  type 2 generated  by an element  $e_{1}$ of weight one  and an
element  $e_{2}$  of  weight  two  subject to  the  relations  $[e_{2}
\U{e_{1}}{i} e_{2}] = 0$, for all  $i \ge 1$.  This is the algebra $m$
of~\eqref{eq:m}.

In positive characteristic, note first of all that in $L$
we may take $e_{1} = z$, $e_{2} = [y z]$ where
$0 \ne y \in C_{2}$, and take $e_{k} = [e_{k-1} e_1]$
for $k > 2$. Suppose that in $M$ we have a segment of the sequence 
of two-step centralizers of the form
\begin{equation*}
  C_{2} = C_{n-2} = C_{n-1}, C_{n} = \Span{y + \lambda z}
  \ne C_{2}, C_{n+1} = C_{n+2} = C_{2},
\end{equation*}
so that $\lambda \ne 0$. 
Note that the first constituent has length $2q \ge 6$ so
that, in particular, 
\begin{equation*}
  [e_{3} e_{2}]
  =
  [ [y z z] [y z] ]
  =
  [y z z y z] - [y z z z y]
  = 0.
\end{equation*}

We have
\begin{align*}
  [e_{n-1} e_{2}] 
  &= 
  [e_{n-1} [y z]]
  \\&=  
  [e_{n-1} y z] - [e_{n-1} z y]
  \\&= 
  - [e_{n-1} z y] && \text{as $C_{n-1} =C_{2} = \Span{y}$}
  \\&= 
  - [e_{n-1} e_{1} y]
  \\&= 
  - [e_{n} y]
  \\&= 
  [e_{n}, \lambda z] - [e_{n}, y + \lambda z]
  \\&=
  \lambda e_{n+1}.
\intertext{Similarly}
  [e_{n} e_{2}] 
  &= 
  [e_{n} [y z]]
  \\&=  
  [e_{n} y z] - [e_{n} z y]
  \\&= 
  [e_{n} y z] && \text{as $C_{n+1} =C_{2} = \Span{y}$}
  \\&= 
  [e_{n} y e_{1}]
  \\&= 
  [e_{n}, - \lambda z, e_{1}] + [e_{n}, y + \lambda z, e_{1}]
  \\&=
  - \lambda e_{n+2}.
\end{align*}
Finally
\begin{equation*}
  [e_{n+1} e_{2}] = 0
\end{equation*}
as $C_{n+1} =  C_{n+2} = C_{2}$.

In view of this, we introduce a definition of constituents for
algebras of type~2 that is compatible with the definition for algebras
of type 1.
Let $L$ be an arbitrary algebra of type 2.
If $[e_{3} e_{2}] = [e_{2} e_{1} e_{2}] \ne 0$, we have
no theory of constituents for $L$. Algebras of this
type are dealt with in Section 3 and Section 7. If
$[e_{2} e_{1} e_{2}] = 0$, and for some $n$ we have
$[e_{n-1} e_{2}] = 0$, but
$[e_{n} e_{2}] = \lambda e_{n+2} \ne 0$, for some
$\lambda \ne 0$, then
\begin{align*}
  0 
  &= 
  [e_{n-1} [e_{2} e_{1} e_{2}]] 
  \\&= 
  - [e_{n-1} e_{1} e_{2} e_{2}]
  + 2 [e_{n-1} e_{2} e_{1} e_{2}]
  - [e_{n-1} e_{2} e_{2} e_{1}]
  \\&=
  - [e_{n-1} e_{1} e_{2} e_{2}]
  \\&=
  - [e_{n} e_{2} e_{2}]
  \\&=
  - \lambda [e_{n+2} e_{2}],
\end{align*}
so that $[e_{n+2} e_{2}] = 0$. We  are therefore led
to the following  definition. Let $L$ be an algebra of
type 1 in which $[e_{2} e_{1} e_{2}] = 0$. Suppose there
are integers $m, n$ such that
\begin{align*}
  &[e_{m-1} e_{2}] = 0,
  \\&
  [e_{m} e_{2}] = \eta e_{m+2}, && \text{with $\eta \ne 0$,}
  \\&
  [e_{m+1} e_{2}] = \theta e_{m+2},
  \\&
  [e_{m+2} e_{2}] = \dots =  [e_{n-1} e_{2}] = 0,
  \\&
  [e_{n} e_{2}] = \lambda e_{n+2}, && \text{with $\lambda \ne 0$,}
  \\&
  [e_{n+1} e_{2}] = \mu e_{n+3}.
\end{align*}
We call this pattern a \emph{constituent of length $l = n - m$
and type $(\lambda,\mu)$.} Note that $\theta$ and $\mu$ might well be
zero.

Here, too, the  first constituent requires an ad hoc
treatment. If in the algebra  $L$ one has
$[e_{2} e_{1}  e_{2}] =  0$, and $n$  is the smallest
integer greater than 1 such that $[e_{n} e_{2}]  \ne 0$, we
say that the first constituent has length $n+1$. If there
is no such $n$, then $L$ is isomorphic to the algebra $m$
above.  

We will see in Section~\ref{sec:first_constituent}
that the first constituent of an algebra of type 2 can
have length $q + 1$ or $2 q$,  where $q$ is a power of the
characteristic of the underlying field. If the first
constituent has length $2 q$, we will see in
Section~\ref{sec:classical} that $L$ comes from an 
algebra of type $1$ via the procedure described above.
If the first constituent has length $q + 1$, we will see in
Sections~\ref{sec:extra}--\ref{sec:construction_3} that we
obtain one soluble algebra of type 2 for $q > 3$, and a family of
soluble algebras for $q = 3$.

We have just seen that an algebra of type 2 that comes
from an algebra of type 1 has constituents of type
$(\lambda,-\lambda)$. We now prove that the converse
also holds.

Suppose all constituents of the algebra $L$ of type 2
are of type $(\lambda,-\lambda)$. Consider the following
partial linear map
\begin{equation*}
  \begin{cases}
    e_{1} \mapsto - e_{2}\\
    e_{2} \mapsto 0.
  \end{cases}
\end{equation*}
We show that we can extend this to a unique derivation $D$ of weight 1
on the whole of  $L$. In the extension $M$ of $L$  by $D$, we have $[D
e_{1}] = - e_{1} D = e_{2}$. Thus $M$ is generated by the elements
$e_{1}$ and $D$ of weight 1, and it is an uncovered algebra of type 1.

We begin with
$e_{3}D=[e_{2} e_{1}]D=[e_{2}D, e_{1}]+[e_{2}, e_{1}D] = 0$.
Suppose now we come to the end of a constituent in $L$,
so that we have
\begin{equation*}
  [e_{i-2} e_{2}] = 0, 
  [e_{i-1} e_{2} ] =  \lambda e_{i+1},
  [e_{i} e_{2} ] = - \lambda e_{i+2}.
\end{equation*}
We have so far, proceeding by induction, $e_{i-2} D = 0$.
Now 
\begin{equation*}
  e_{i-1} D = [e_{i-2} e_{1}] D = 
  [e_{i-2} D, e_{1}] + [e_{i-2}, e_{1} D]
  = - [e_{i-2} e_{2}] = 0.
\end{equation*}
Then
\begin{equation*}
  e_{i} D = [e_{i-1} e_{1}] D = 
  [e_{i-1} D, e_{1}] + [e_{i-1}, e_{1} D]
  = - [e_{i-1} e_{2}] = - \lambda e_{i+1},
\end{equation*}
\begin{equation*}
  e_{i+1} D = [e_{i} e_{1}] D = 
  [e_{i} D, e_{1}] + [e_{i}, e_{1} D]
  = - \lambda [e_{i+1} e_{1}] - [e_{i} e_{2}] = 0,
\end{equation*}
and
\begin{equation*}
  e_{i+2} D = [e_{i+1} e_{1}] D = 
  [e_{i+1} D, e_{1}] + [e_{i+1}, e_{1} D]
  = - [e_{i+1} e_{2}] = 0,
\end{equation*}
so that we can continue by induction.

This definition of $D$ is compatible with
the relations $[e_{i-2}, e_{2}] = 0$,
$[e_{i-1}, e_{2}] = - \lambda e_{i+1}$,
$[e_{i}, e_{2}] = \lambda e_{i+2}$,
$[e_{i+1}, e_{2}] = 0$. 
This is clear for all but the third one. For this we have
\begin{align*}
  [e_{i} D, e_{2}] + [e_{i}, e_{2} D] = - \lambda [e_{i+1}, e_{2}]
  = 0 = e_{i+2} D.
\end{align*}

In~\cite{CMN} a device for studying algebras of
type 1 called \emph{deflation} has been introduced.
We now show that this can be applied also to algebras
of type 2, and the result will be an algebra of
type 1. This is useful in simplifying some proofs later on.

Let $L$ be an algebra of type 2 as in~\eqref{eq:L_again}.
Consider its subalgebra
\begin{equation*}
  S = \bigoplus_{i = 1}^{\infty} L_{i p}.
\end{equation*}
Grade $S$ by assigning weight $i$ to $L_{i  p}$.
Now $S$ admits the derivation $D = \ad(e_{1})^{p}$ which,
in the new grading, has weight 1. We have 
\begin{equation*}
  L_{i p} \ad(e_{1})^{p} = [L_{i p} \U{e_{1}}{p}] = L_{(i+1) p}.
\end{equation*}
It  follows that the extension of $S$ by  $D$ is a graded
Lie algebra of maximal class, and it is generated by the
two elements $e_{p}$ and $D$ of weight 1. Therefore it is
an algebra of type 1.

\bigskip

In this section we have used several times the Jacobi identity $[z [y
x]] = [z y x] - [z x y]$, and its consequence
\begin{equation*}
  [z [y \U{x}{n}]] = \sum_{i = 0}^{n}
  (-1)^{i} \binom{n}{i} [z \U{x}{i} y \U{x}{n-i}].
\end{equation*}
In such a formula, to evaluate binomial coefficients modulo a prime we
will make use of Lucas' theorem, in the following form. Suppose $a, b$
are non-negative  integers, and $q  > 1$  is a power  of a prime  $p$. 
Write $a  = a_{0} +  a_{1} q$, and  $b = b_{0}  + b_{1} q$,  where the
$a_{i}$ and $b_{i}$ are non-negative integers, and $a_{0}, b_{0} < q$.
Then
\begin{equation*}
  \binom{a}{b} \equiv \binom{a_{0}}{b_{0}} \cdot \binom{a_{1}}{b_{1}}
  \pmod{p}.
\end{equation*}

\section{Characterizing $m_{2}$}

In this section we start dealing with algebras
of type 2 that do not admit a theory of constituents,
that is, in which $[e_3,e_2] \ne 0$ We may
thus assume without loss of generality $[e_3,e_2]=e_5$.
We obtain 
\begin{align*}
  0
  &=
  [e_3 e_3]
  \\&=
  [e_3 [e_2 e_1]]
  \\&=
  [e_3 e_2 e_1] - [e_3 e_1 e_2]
  \\&=
  e_6 - [e_4 e_2].
\end{align*}

Suppose that
\[
\lbrack e_i,e_1]=e_{i+1}\text{ for }i>1, 
\]
\[
\lbrack e_3,e_2]=e_5,\;[e_4,e_2]=e_6,\;[e_5,e_2]=ae_7,%
\;[e_6,e_2]=be_8, 
\]
\[
\lbrack
e_7,e_2]=ce_9,\;[e_8,e_2]=de_{10},\;[e_9,e_2]=fe_{11},%
\;[e_{10},e_2]=ge_{12}. 
\]
Here $a,b,c,d,f,g,$ are parameters.

$[[e_2,e_1,e_1],[e_2,e_1,e_1]]=0$ gives $1-2a+b=0$,
so $b=-1+2a$.

$[[e_2,e_1,e_1,e_1],[e_2,e_1,e_1,e_1]]=0$ gives
$a-3b+3c-d=0$.

Now $a-3b+3c-d=-5a+3+3c-d$, and so $d=3-5a+3c$.

$[[e_2,e_1,e_1,e_1,e_1],[e_2,e_1,e_1,e_1,e_1]]=0$
gives $b-4c+6d-4f+g=0$.
\[
b-4c+6d-4f+g=17-28a+14c-4f+g, 
\]
so $g=-17+28a-14c+4f$.

Note that $[e_2,e_1,e_1,e_1]=-[e_1,e_2,e_2]$.

$[[e_2,e_1,e_1,e_1],[e_1,e_2,e_2]]=0$ gives $bd-2ad+ac=0$. 
\[
\begin{array}{l}
bd-2ad+ac \\ 
=\left(-1+2a\right) \left(3-5a+3c\right)-2a\left(3-5a+3c\right)+ac \\ 
=-3+5a-3c+ac,
\end{array}
\]

so either $a=3$ (which gives $12=0$), or
$\displaystyle c=\frac{3-5a}{a-3}$.

\[
\lbrack
[e_2,e_1,e_1,e_1,e_1],[e_2,e_1,e_1,e_1]]+[[e_2,e_1,e_1,e_1,e_1],[e_1,e_2,e_2]]=0 
\]
gives 
\[
b-3c+3d-f+cf-2bf+bd=0. 
\]
\begin{align*}
  \displaystyle 
  b-3c+3d-f+cf&-2bf+bd =
  \\&=
  \displaystyle 
  8-13a+6\frac{3-5a}{a-3}-f+\frac{3-5a}{a-3}f
  \\&\phantom{=\ }
  -2\left( -1+2a\right) f
  +
  \allowbreak \left( -1+2a\right) \left( 3-5a+3\frac{3-5a}{a-3}\right)
  \\&=
  \displaystyle 
  -2\frac{-7a+3+a^2-4fa+2fa^2+5a^3}{a-3}.
\end{align*}
So provided the characteristic is not 2, and provided $a\neq 3$, 
\[
-7a+3+a^2+5a^3+(2a^2-4a)f=0. 
\]

\[
\lbrack
[e_2,e_1,e_1,e_1,e_1,e_1],[e_2,e_1,e_1,e_1]]+[[e_2,e_1,e_1,e_1,e_1,e_1],[e_1,e_2,e_2]]=0 
\]
gives 
\[
c-3d+3f-g+dg-2cg+cf=0. 
\]
\[
\begin{array}{l}
\displaystyle 
c-3d+3f-g+dg-2cg+cf \\ 
\displaystyle 
=-2 \frac{18-27f-123a+280a^2+78fa-53fa^2-253a^3+10fa^3+70a^4}{\left(
a-3\right) ^2}.
\end{array}
\]
So provided the characteristic is not 2, and provided $a\neq 3$, 
\[
18-123a+280a^2-253a^3+70a^4+(-27+78a-53a^2+10a^3)f=0. 
\]
Combining these two equations we obtain 
\begin{multline*}
(-27+78a-53a^2+10a^3)(-7a+3+a^2+5a^3)
-\\-
(2a^2-4a)(18-123a+280a^2-253a^3+70a^4)=0 
\end{multline*}
Expanding, we obtain 
\begin{align*}
  0
  &=
  495a-81-1260a^2+1710a^3-1305a^4+531a^5-\allowbreak 90a^6
  \\&=
  -9\times \left(10a-9\right) \left( a-1\right) ^5. 
\end{align*}
So if the characteristic is not 2 or 3 or 5 then
$a=1$ or $a=\frac 9{10}$. If the characteristic is 5
then $a=1$. The cases when the characteristic is
2 or 3 have to be dealt with separately. We deal
with the latter in Section~\ref{sec:extra_3}. 

When $a=\frac 9{10}$, it is proved in~\cite{Marina}
that the algebras one obtains are quotients of a
certain central extension of the positive part of
the infinite-dimensional Witt algebra. In any case,
there are no infinite-dimensional algebras of maximal
class here.

The choice $a = 1$ uniquely determines the following
metabelian Lie algebra \cite{F, SZ}:
\begin{align*}
  m_{2} = 
  \big\langle 
  e_{i}, i \ge 1 :\ 
  &\text{$[e_{i} e_{1}] = e_{i+1}$, for $i \ge 2$,}
  \\&\text{$[e_{i} e_{2}] = e_{i+2}$, for  $i \ge 3$,}
  \\&\text{$[e_{i} e_{j}] = 0$, for $i, j \ge 3$ }
  \big\rangle.
\end{align*}
Note that $\ad(e_{2})$ is the square of $\ad(e_{1})$
on $L^{2}$.

In fact, we have to show that $m_{2}$ has the
following presentation:
\begin{equation*}
  \Span{e_{1}, e_{2} : [e_{2} e_{1} e_{2}] = [e_{2} e_{1}^{3}],
    [e_{2} e_{1}^{3} e_{2}] = [e_{2} e_{1}^{5}]}.
\end{equation*}
We use the notation $[e_{i} e_{1}] = e_{i+1}$,
so that the two defining relations can be rewritten
as $[e_{3} e_{2}] = e_{5}$ and $[e_{5} e_{2}] = e_{7}$.
We have already seen that the first one implies
$[e_{4} e_{2}] = e_{6}$. Suppose now we have proved
\begin{equation*}
  [e_{3} e_{2}] = e_{5},
  [e_{4} e_{2}] = e_{6},
  [e_{5} e_{2}] = e_{7},
  \dots,
  [e_{n-1} e_{2}] = e_{n+1},
\end{equation*}
for some $n > 5$, and want to prove
$[e_{n} e_{2}] = e_{n+2}$. We work out the expansion
\begin{align*}
  0 
  &=
  [ e_{n-3}, [e_{3} e_{2}] - e_{5}]
  =
  [ e_{n-3}, [e_{2} e_{1} e_{2}] - [e_{2} e_{1}^{3}]]
  \\&=
  [ e_{n-3} [e_{2} e_{1}] e_{2}]
  \\&\phantom{=\ }
  -  [ e_{n-3}  e_{2} [e_{2} e_{1}]]
  \\&\phantom{=\ }
  -  [ e_{n-3} [e_{2} e_{1}^{3}]
  \\&=
  (1 - 1) [e_{n} e_{2}]
  \\&\phantom{=\ }
  - e_{n+2} + [e_{n} e_{2}]
  \\&\phantom{=\ }
  - (1 - 3 + 3) e_{n+2} + [e_{n} e_{2}]
  \\&=
  2 [e_{n} e_{2}] - 2 e_{n+2}.
\end{align*}
Note that this does not work for $n = 5$. From this it is
straightforward to see that the algebra is metabelian,
and thus is isomorphic to $m_{2}$. In fact we have for
$i, j \ge 3$
\begin{align*}
  [ e_{i} e_{j} ] 
  &=
  [ e_{i} [ e_{2} \U{e_{1}}{j-2} ] ]
  \\&=
  \sum_{k=0}^{j-2} (-1)^{k} \binom{j-2}{k} 
    [ e_{i} \U{e_{1}}{k} e_{2} \U{e_{1}}{j-2-k} ]
  \\&=
  \left(  
    \sum_{k=0}^{j-2} (-1)^{k} \binom{j-2}{k}
  \right)
  \cdot
  e_{i+j}
  \\&=
  0.
\end{align*}

\section{The length of the first constituent}
\label{sec:first_constituent}

Suppose  now $L$  is an  algebra of  type 2 over a field of positive
characteristic $p$. Suppose $L$ admits  a  theory of
constituents. Therefore $[e_{3} e_{2}] = [e_{2} e_{1} e_{2}] = 0$. If
$[e_{i}  e_{2}] = 0$  for all  $i \ge  3$, then  $L$ is  isomorphic to
$m$ of~\eqref{eq:m}. Suppose thus there is an $n > 3$ such that
\begin{math}
  [e_{3} e_{2}] = [e_{4} e_{2}] = \dots = [e_{n-2} e_{2}] = 0,
\end{math}
but
\begin{math}
  [e_{n-1} e_{2}] \ne 0.
\end{math}
We intend to  show that $n$, the length of  the first constituent, can
only assume the values
\begin{equation*}
\begin{cases}
  2 q, & \text{for some power $q$ of $p$, or}\\
  q + 1, & \text{for some power $q > 3$ of $p$.}
\end{cases}
\end{equation*}

We may assume, rescaling $e_{2}$, that $[e_{n-1} e_{2}] = e_{n+1}$. We
first prove  that $n$ is even,  with a simple argument  similar to one
of~\cite{CMN}. In fact, if $n = 2 k - 1$ is odd, we have
\begin{align*}
  0
  &=
  [ e_{k} e_{k} ] = [e_{k} [e_{2} \U{e_{1}}{k-2}]]
  \\&=
  \sum_{i = 0}^{k-2} 
  (-1)^{i} \binom{k-2}{i}
  [ e_{k} \U{e_{1}}{i} e_{2} \U{e_{1}}{k-2-i} ]
  \\&=
  \sum_{i = 0}^{k-2} 
  (-1)^{i} \binom{k-2}{i}
  [ e_{k+i} e_{2} \U{e_{1}}{k-2-i} ]
  \uptoasign
  [e_{n-1} e_{2}] 
  = 
  e_{n+1},
\end{align*}
a contradiction. 
Here and in  the following we write $a \uptoasign b$  to mean that $a$
is either $b$ or  $- b$.  Write $n = 2 k$. We  aim at proving that the
only possible values for $k$ are $q$ and $(q+1)/2$.

We first compute
\begin{align*}
  0
  &=[e_{k+1} e_{k+1} ] =
  [e_{k+1} [e_{2} \U{e_{1}}{k-1}] ]
  \\&\uptoasign
  \binom{k - 1}{k - 2} [e_{2 k - 1} e_{2} e_{1}]
  - \binom{k - 1}{k -1} [e_{2 k} e_{2}]
  \\&=
  (k - 1) e_{2 k + 2} - [e_{2 k} e_{2}],
\end{align*}
to show
\begin{equation*}
  [e_{n} e_{2}]
  =
  (k - 1) e_{n + 2}.
\end{equation*}

We now have
\begin{equation}\label{eq:extend_by_one}
  \begin{aligned}
    0
    &=
    [e_{n - 2} [e_{2} e_{1} e_{2}]]
    \\&\uptoasign
    [e_{n - 2} e_{1} e_{2} e_{2}]
    \\&=
    [e_{n + 1} e_{2}].
  \end{aligned}
\end{equation}
Further,
\begin{equation}\label{eq:previous_calc}
  \begin{aligned}
    0
    &=
    [e_{n - 1} [e_{2} e_{1} e_{2}]]
    \\&\uptoasign
    [e_{n - 1} e_{1} e_{2} e_{2}]
    - 2   [e_{n - 1} e_{2} e_{1} e_{2}]
    + [e_{n - 1} e_{2} e_{2} e_{1}]
    \\&=
    (k - 1 - 2) [e_{n+2} e_{2}].
  \end{aligned}
\end{equation}
This  shows that  $[e_{n+2}  e_{2}] =  0$,  except when  $k \equiv  3
\pmod{p}$. 

Suppose first we have $n = 6$, or $k = 3$. We have here 
\begin{equation*}
  [e_{5} e_{2}] = e_{7},\quad
  [e_{6} e_{2}] = 2 e_{8}, \quad
  [e_{7} e_{2}] = 0.
\end{equation*}
We want to show  that $p = 3$ or $5$ here, so  that this fits into the
$n = 2q$ or $q+1$ pattern above. Suppose $p > 5$. We compute
\begin{align*}
  0 
  &= 
  [e_{5} [e_{2} e_{1}^{3}] ]
  \\&=
  [e_{5} e_{2} e_{1}^{3} ] - 3   [e_{6} e_{2} e_{1}^{2} ]
  + 3  [e_{7} e_{2} e_{1} ] -   [e_{8} e_{2} ]
  \\&=
  -5 e_{10}  - [e_{8} e_{2}],
\end{align*}
so that $[e_{8} e_{2}] = -5 e_{10}$.
\begin{align*}
  0 
  &= 
  [e_{6} [e_{1} e_{2} e_{2}] ]
  \\&=
  [e_{7} e_{2} e_{2} ] - 2   [e_{6} e_{2} e_{1} e_{2} ]
  + [e_{6} e_{2}  e_{2} e_{1} ]
  \\&=
  - 4 [e_{9} e_{2}]  + 2 (-5) e_{11},
\end{align*}
so that
\begin{equation*}
  [e_{9} e_{2}] = - \frac{5}{2} e_{11}.
\end{equation*}
Finally
\begin{align*}
  0 
  &= 
  [e_{6} [e_{2} e_{1}^{4}] ]
  \\&=
  [e_{6} e_{2} e_{1}^{4} ]
  - 4  [e_{7} e_{2} e_{1}^{3} ] + 6 [e_{8} e_{2} e_{1}^{2} ]
  - 4 [ e_{9} e_{2} e_{1}] + [e_{10} e_{2} ]
  \\&=
  18 e_{12}  + [e_{10} e_{2}]
\end{align*}
and
\begin{align*}
  0 
  &= 
  [e_{7} [e_{1} e_{2} e_{2}] ]
  \\&=
  [e_{8} e_{2} e_{2} ] - 2   [e_{7} e_{2} e_{1} e_{2} ]
  + [e_{7} e_{2}  e_{2} e_{1} ]
  \\&=
  -5 [e_{10} e_{2}]
\end{align*}
yield $e_{12} = 0$, a contradiction.

Suppose then $k > 3$, that is, $n > 6$. We have thus $[e_{5} e_{2}] =
0$, so that 
\begin{align*}
  0
  &=
  [e_{n - 3} [e_{5} e_{2}]]
  =
  [e_{n - 3} [e_{2} e_{1} e_{1} e_{1} e_{2}] ]
  \\&=
  [e_{n - 3} [e_{2} e_{1} e_{1} e_{1}] e_{2}]]
  \\&=
  3 [e_{n - 3} e_{1} e_{1} e_{2} e_{1} e_{2} ]
  -   [e_{n - 3} e_{1} e_{1} e_{1} e_{2} e_{2}]
  \\&=
  (3 - (k-1)) [e_{n+2} e_{2}].
\end{align*}
This shows that $[e_{n+2}  e_{2}] =  0$,  except when $k \equiv 4
\pmod{p}$, which was covered by~\eqref{eq:previous_calc}.

To find out what the possible values of $k$ are, we compute
\begin{align*}
  0
  &=
  [e_{k+2} e_{k+2}] = [e_{k+2} [e_{2} \U{e_{1}}{k}] ]
  \\&\uptoasign
  \binom{k}{k - 3} [e_{2 k - 1} e_{2} e_{1} e_{1} e_{1}]
  - \binom{k}{k - 2} [e_{2 k} e_{2} e_{1} e_{1}]
\end{align*}
which yields
\begin{equation*}
 0 
 =
 \left( \binom{k}{3}  - \binom{k}{2} (k-1) \right) e_{2k+4}
 = \frac{k (k - 1) (-2k + 1)}{6} \, e_{2k+4}.
\end{equation*}
This shows  that the only  possibilities for $k$  are 
\begin{equation*}
  k \in \left\{  0, 1,\frac{1}{2} \right\} \pmod{p},
\end{equation*}
for $p > 3$, whereas  for $p = 3$ one has
\begin{equation*}
  k \in \left\{  0, 1, \frac{1}{2} \right\} \pmod{9}.
\end{equation*}

When $k \equiv 0 \pmod{p}$, we show that $k = q$, a power of $p$. (The
case when $p =  3$ is not special here, as we  have already dealt with
$k =  3$ for $p =  3$ above.) This  we do by exploiting  the deflation
procedure,  as described in  Section~\ref{sec:preliminaries}.  Suppose
in  fact $k =  q m$,  with $q$  a power  of $p$,  and $m  \not\equiv 0
\pmod{p}$.  Thus $n = 2 q m$ here. We  have $[e_{n - 1} e_{2}] = e_{n + 1}$
and $[e_{n} e_{2}] = - e_{n + 2}$.  We have also proved
in~\eqref{eq:extend_by_one} that $[e_{n+1} e_{2}]  = 0$.
We first extend this to
\begin{equation*}
  [e_{n + 1} e_{2}] 
  = 
  [e_{n + 2} e_{2}]
  =
  \dots
  =
  [e_{n + p - 2} e_{2}] 
  =
  0.
\end{equation*}
We proceed by induction on $l$, for $1 < l \le p - 2$:
\begin{equation}\label{eq:implicit_lower_bound}
  \begin{aligned}
    0
    &=
    [e_{n - 1} [e_{2} e_{1}^{l - 1} e_{2}]]
    \\&=
    [e_{n - 1} [e_{2} e_{1}^{l - 1}] e_{2}]
    -
    [e_{n - 1} e_{2} [e_{2} e_{1}^{l - 1}]]
    \\&=
    [e_{n - 1} e_{2} e_{1}^{l - 1} e_{2}]
    - (l - 1)
    [e_{n - 1} e_{1} e_{2} e_{1}^{l-2} e_{2}]
    - (-1)^{l-1}
    [e_{n - 1} e_{2}  e_{1}^{l - 1} e_{2}]
    \\&=
    (1 + l - 1 - (-1)^{l - 1})
    \cdot
    [e_{n + l} e_{2}].
  \end{aligned}
\end{equation}
Now
\begin{equation*}
  1 + l - 1 - (-1)^{l - 1} = l + (-1)^{l} =
  \begin{cases}
    l - 1 & \text{when $l$ is odd,}\\
    l + 1 & \text{when $l$ is even.}
  \end{cases}
\end{equation*}
In any case the coefficient of $ [e_{n + l} e_{2}]$ is less than $p$
for $l < p - 1$, so that it is non-zero.

In the deflated algebra, we thus have
\begin{equation*}
  [e_{2 q m - p} e_{p}] =
  [e_{2 q m - p} [e_{2} e_{1}^{p-2}]] =
  0
\end{equation*}
and
\begin{equation*}
  [e_{2 q m} e_{p}] =
  [e_{2 q m} [e_{2} e_{1}^{p-2}]] =
  - e_{2 q m + p}.
\end{equation*}
In the deflated algebra the first constituent has thus length $2 q m$.
It follows from the theory of algebras of type~1 that $m  = 1$.  We will
show in Section \ref{sec:classical} that  algebras of type 2 with $k =
q$ come from algebras of type 1.

When $k \equiv \displaystyle \frac{1}{2} \pmod{p}$, write
\begin{equation*}
  k = \frac{q m + 1}{2},
\end{equation*}
where $p$  does not divide $m$.   Thus $n = q  m + 1$. We  want to show
that $m = 1$. Suppose otherwise.  We have
\begin{equation*}
  [e_{n - 1} e_{2}] = e_{n + 1}
  \qquad\text{and}\qquad
  [e_{n}  e_{2}] = - \frac{1}{2} e_{n + 2}.
\end{equation*}  
We
begin with proving 
\begin{equation}\label{eq:plenty_of_zeroes}
  [e_{n + 1} e_{2}] 
  = 
  [e_{n + 2} e_{2}] 
  = \dots = 
  [e_{n +  q - 1} e_{2}] 
  =  
  0.
\end{equation}

The identity
\begin{equation}\label{eq:an_equation}
  [e_{2} e_{1}^{k} e_{2}] = 0
\end{equation}
holds for $k \le n - 4$. Note that $n - 4 = q  m - 3 \ge 2 q - 3$,
as $m > 1$.

Let $l < q - 1$. Write $l + 1 = \beta p^{t}$, where $\beta \not\equiv 0
\pmod{p}$. Note that $p^{t} < q$, so that
\begin{equation*}
  l + p^{t} \le q - 2 + p^{t} \le 2 q - 3, 
\end{equation*}
and $[e_{2} e_{1}^{l + p^{t}} e_{2}] = 0$, by~\eqref{eq:an_equation}.

Suppose first $t > 0$. 
We compute
\begin{align*}
  0
  &=
  [ e_{n - 1 - p^{t}} [e_{2} e_{1}^{l + p^{t}} e_{2}] ]
  \\&=
  [ e_{n - 1 - p^{t}} [e_{2} e_{1}^{l + p^{t}}]  e_{2}] 
  \\&\uptoasign
  \left(
    \binom{l + p^{t}}{p^{t}} 
    + 
    \frac{1}{2}
    \binom{l + p^{t}}{p^{t} + 1} 
  \right)
  [e_{n + l + 1} e_{2}]
  \\&=
  \left(
    \binom{\beta p^{t} + p^{t} - 1}{p^{t}} 
    + 
    \frac{1}{2}
    \binom{\beta p^{t} + p^{t} - 1}{p^{t} + 1} 
  \right)
  [e_{n + l + 1} e_{2}]
  \\&=
  \left( \beta + \frac{1}{2} (\beta \cdot (-1)) \right)
  [e_{n + l + 1} e_{2}]
  \\&=
  \frac{\beta}{2}
  \cdot
  [e_{n + l + 1} e_{2}],
\end{align*}
so that $[e_{n + l + 1} e_{2}] = 0$.

Now consider the  case when $p^{t} =  1$, so that $l +  1 \not\equiv 0
\pmod{p}$. An analogous calculation yields
\begin{equation*}
  0
  =
    \frac{(l + 1) \cdot (l + 4)}{4}
    \cdot  
    [e_{n + l + 1} e_{2}].
\end{equation*}
We obtain $[e_{n + l + 1} e_{2}] = 0$, except when $l + 4$ is divisible
by $p$. Note that we may assume $p > 3$ here, since we have already
dealt with the case when $l + 1 \equiv 0 \pmod{p}$. We compute
\begin{align*}
  0
  &=
  [e_{n-3} [e_{2} e_{1}^{l + 2} e_{2}]]
  \\&\uptoasign
  \left(
    \binom{l + 2}{2} + \frac{1}{2} \cdot \binom{l + 2}{3}
  \right)
  [e_{n + l + 1} e_{2}]
  \\&=
  [e_{n + l + 1} e_{2}],
\end{align*}
as $p > 3$.

We now reach a  contradiction by proving $e_{n + q +1} =  e_{q m + q +
  2} = 0$. Since
$n = q m + 1$ is even, $m$ is odd, and $q m + q + 2$ is even. Consider
the integer
\begin{equation*}
  \frac{q m + q + 2}{2} 
  =
  q \cdot \frac{m+1}{2} + 1.
\end{equation*}
Note that 
\begin{equation*}
   q \cdot \frac{m+1}{2} + 1 - 2 = 
   q \cdot  \frac{m-1}{2} + q - 1
\end{equation*}
We obtain, using~\eqref{eq:plenty_of_zeroes},
\begin{equation}\label{eq:the_big_kill}
  \begin{aligned}
    0
    &=
    [ 
    e_{q \cdot \frac{m+1}{2} + 1}
    e_{q \cdot \frac{m+1}{2} + 1}
    ]
    =
    [
    e_{q \cdot \frac{m+1}{2} + 1}
    [e_{2}e_{1}^{q \frac{m-1}{2} + q - 1}]
    ]
    \\&=
    \begin{aligned}[t]
      \Big(
      &(-1)^{q \frac{m-1}{2} - 1} 
      \binom{q \frac{m-1}{2} + q - 1}{q \frac{m-1}{2} - 1}
      \\&+
      (-1)^{q \frac{m-1}{2}} 
      \binom{q \frac{m-1}{2} + q - 1}{q \frac{m-1}{2}} 
      \cdot \left(  - \frac{1}{2}\right)
      \Big)
      \cdot
      e_{q m + q + 2}.
    \end{aligned}
  \end{aligned}
\end{equation}
Now we have
\begin{equation*}
  \binom{q \frac{m-1}{2} + q - 1}{q \frac{m-1}{2} - 1}
  \equiv
  \binom{q \frac{m-1}{2} + q - 1}{q \frac{m-3}{2} + q - 1}
  \equiv
  \frac{m-1}{2}
  \pmod{p},
\end{equation*}
while
\begin{equation*}
  \binom{q \frac{m-1}{2} + q - 1}{q \frac{m-1}{2}}
  =
  1.
\end{equation*}
Therefore, up to a sign, the overall coefficient of $e_{q m + q + 2}$
in~\eqref{eq:the_big_kill} is
\begin{equation*}
  \frac{m-1}{2} + \frac{1}{2} = \frac{m}{2} \not\equiv 0.
\end{equation*}
This disposes of the case $m > 1$, so we obtain
\begin{equation*}
  k = \frac{q + 1}{2},
  \quad
  n = q + 1.
\end{equation*}
We will deal with this case
in Sections \ref{sec:extra} and \ref{sec:Extra_Construction}.
Remember that when $p = 3$ we are taking $q \ge 9$ here. In fact when
$q =3$ we get $k = 2$,  so that $[e_{3} e_{2}] \ne 0$, and the algebra
does not admit a theory of  constituents.

We now  deal with the case  $k \equiv 1 \pmod{p}$,  so $k = 1  + q m$,
where $q$ is a power of $p$,  and $m \not\equiv 0 \pmod{p}$. Thus $n =
2 q m +  2$. We have thus $[e_{n - 1} e_{2}] =  e_{n + 1}$ and $[e_{n}
e_{2}] = 0$. We want to show that this case does not occur.

Let $1 \le l < q$. Assume by induction
\begin{equation*}
  [e_{n} e_{2}] = [e_{n+1} e_{2}] = \dots = [e_{n+l-1} e_{2}] = 0.
\end{equation*}
We compute
\begin{equation}\label{eq:l}
  0 
  = 
  [e_{n-2} [e_{2} e_{1}^{l} e_{2}]]
  =
  [e_{n-2} [e_{2} e_{1}^{l}] e_{2}]
  =
  - l [e_{n + l} e_{2}].
\end{equation}
We obtain $[e_{n  + l} e_{2}] = 0$  for $l < p$.  We can  use this and
deflation to show that $m = 1$. 
Because of
\begin{equation*}
  [e_{n - 2} e_{p}] 
  =
  [e_{n - 2} [e_{2} e_{1}^{p-2}]]
  =
  2 e_{n + p - 2},
\end{equation*}
the length of the first constituent in the deflated algebra (which is
of type~1) is $2 q m / p$. If $m > 1$, this is not twice a power of
$p$. It follows that $m = 1$, and $n = 2 q + 2$.

We now show  that $[e_{n + l} e_{2}]  = 0$ holds in fact for  all $l <
q$. Because of the argument  of~\eqref{eq:l}, we have to deal with the
case $l  \equiv 0 \pmod{p}$.  If $p^{t}$ is  the highest power  of $p$
that  divides $l$, and  $l =  \beta p^{t}$,  with $\beta  \not\equiv 0
\pmod{p}$, we compute
\begin{align*}
  0
  &=
  [e_{n-p^{t}-1} [e_{2} e_{1}^{l+p^{t}-1} e_{2}]]
  \\&=
  [e_{n-p^{t}-1} [e_{2} e_{1}^{l+p^{t}-1}] e_{2}]
  \\&=
   \pm \binom{l + p^{t} - 1}{p^{t}} [e_{n + l} e_{2}].
\end{align*}
Here 
\begin{equation*}
  \binom{l + p^{t} - 1}{p^{t}} 
  = 
  \binom{\beta p^{t} + p^{t} - 1}{p^{t}} \equiv \beta
  \not\equiv 0 \pmod{p}.
\end{equation*}
We can perform this calculation when $l +  r - 1 < 2 q - 1$. Note that
this holds for $l < q$. We have thus proved
\begin{equation}\label{eq:these_vanish}
  [e_{n} e_{2}] = [e_{n + 1} e_{2}] = 
  \dots = [e_{n + q - 1} e_{2}] = 0.
\end{equation}
Now we use the  relation $[e_{n - 1} e_{2}] - e_{n +  1} = 0$ to prove
$e_{3 q + 3} = e_{n + q + 1} = 0$, a contradiction. We evaluate
\begin{align*}
  0
  &=
  [e_{q}, [e_{n-1} e_{2}] - e_{n+1}]
  \\&=
  [e_{q} [e_{2} e_{1}^{2 q - 1} e_{2}]] - [e_{q} [e_{2} e_{1}^{2 q +
    1}]]
\end{align*}
Note first that $2 q + 1$ is the only value $i$ in the range $2 \le  i
\le 3  q + 1$  for which $[e_{i}  e_{2}] \ne 0$. Now $[e_{q} [e_{2}
e_{1}^{2 q - 1} e_{2}]]$ expands as a combination of commutators of the
form $[e_{i} e_{2} e_{1}^{2 q + 1 - i}]$, for some $q + 2 \le i \le 3
q + 1$, so that it vanishes.
We obtain
\begin{align*}
  0
  &=
  [e_{q} [e_{2} e_{1}^{2 q + 1}]]
  \\&=
  (-1)^{q+1} \binom{2 q + 1}{q+1}
  e_{3 q + 3} 
  \\&=
  2 e_{3 q + 3}.
\end{align*}

\section{First constituent of length $2 q$}
\label{sec:classical}

This is the case $k = q$ of the previous section. Suppose we have
\begin{gather*}
   [e_{i} e_{2}] = 0, \quad  \text{for $i < 2 q - 1$}\\
   [e_{2 q - 1} e_{2}] = e_{2 q + 1},
   \qquad
   [e_{2 q} e_{2}] = - e_{2 q + 2}.
\end{gather*}
We want to show that the algebra comes from an algebra of type $1$ via
the procedure described in Section~\ref{sec:preliminaries}, by proving
that all constituents have type $(\lambda, - \lambda)$.  

Proceeding by
induction, assume  we have  already  proved  this  up to  a  certain
constituent, that ends as
\begin{equation}\label{eq:end_of_constituent}
  [e_{m} e_{2}] = \lambda e_{m+2},
  \qquad
  [e_{m+1} e_{2}] = - \lambda e_{m+3},
\end{equation}
for some $\lambda  \ne 0$.  We first show, also  by induction, that $2
q$ is an  upper bound for the length of the  next constituent, and $q$
is a lower bound.

Suppose the next constituent has length greater than $2 q$, so that
\begin{equation*}
  [e_{m+k} e_{2}] = 0
\end{equation*}
for $2 \le k \le 2 q$. We obtain immediately
\begin{equation*}
  [e_{m} [e_{2 q} e_{2}]]
  =
  [e_{m} [e_{2} e_{1}^{2 q - 2} e_{2}]] 
  = 
  0,
\end{equation*}
as this is a multiple of $[e_{m + 2 q} e_{2}] = 0$.
This yields
\begin{align*}
  0
  &=
  [ e_{m}, [e_{2 q} e_{2}] + e_{2 q + 2}]
  \\&=
  [ e_{m} e_{2 q + 2}]
  =
  [ e_{m} [e_{2} e_{1}^{2 q}]]
  =
  [ e_{m} e_{2} e_{1}^{2 q}]
  \\&=
  \lambda e_{m + 2 + 2 q},
\end{align*}
a contradiction.

We now prove that  the next constituent has length at  least $q$, that
is,
\begin{equation*}
  [e_{m + 2} e_{2}] = [e_{m + 3} e_{2}] = \dots
  = [e_{m + q - 1} e_{2}] = 0.
\end{equation*}
This we do more generally for the case when the current constituent is
of the general form
\begin{equation}\label{eq:general_end_of_constituent}
  [e_{m} e_{2}] = \mu e_{m+2},
  \qquad
  [e_{m+1} e_{2}] = \nu e_{m+3},
\end{equation}
as this will be useful later in this section. Recall that $\mu \ne 0$
here, but $\nu$ might be zero.

If  $\nu = 0$  in ~\eqref{eq:general_end_of_constituent},  we compute,
proceeding by induction on $l$, for $0 < l < q - 1$,
\begin{align*}
  0
  &=
  [e_{m-1} [e_{2} e_{1}^{l} e_{2}]]
  \\&=
  [e_{m-1} [e_{2} e_{1}^{l}] e_{2}]
  \\&=
  - l \mu [e_{m + l + 1} e_{2}]. 
\end{align*}
The coefficient  vanishes when  $l \equiv 0  \pmod{p}$. In  this case,
write $l = \beta p^{t}$, with $\beta \not\equiv 0 \pmod{p}$. Note that
$p^{t} < q$ here,  so that $l + p^{t} - 1 <  q - 2 + q - 1  < 2 q - 3$
and $[e_{2} e_{1}^{l+p^{t}-1} e_{2}] = 0$. Also, $[e_{m-p^{t}} e_{2}]
= \dots = [e_{m-1} e_{2}] = 0$, since we are assuming by induction
that constituents have length at least $q$. We compute
\begin{align*}
  0
  &=
  [e_{m-p^{t}} [e_{2} e_{1}^{l+p^{t}-1} e_{2}]]
  \\&
  - \binom{l + p^{t} - 1}{p^{t}} \mu [e_{m + l + 1} e_{2}]. 
\end{align*}
Here 
\begin{equation*}
  \binom{l + p^{t} - 1}{p^{t}} 
  =
  \binom{\beta p^{t} + p^{t} - 1}{p^{t}}
  \equiv
  \beta
  \not\equiv
  0
  \pmod{p}.
\end{equation*}

Suppose now $\nu \ne 0$. We have first
\begin{align*}
  0
  =
  [e_{m-1} [e_{2} e_{1} e_{2}]]
  =
  - \mu [e_{m+2} e_{2}],
\end{align*}
so that  $[e_{m+2} e_{2}] = 0$.   We proceed now by  induction on $l$,
for $0 < l < q - 2$.
\begin{align}
  0
  &=
  [e_{m} [e_{2} e_{1}^{l} e_{2}]]
  \notag
  \\&=
  [e_{m} [e_{2} e_{1}^{l}] e_{2}]
  -
  [e_{m} e_{2} [e_{2} e_{1}^{l-1}]]
  \notag
  \\&=
  (\mu - l \nu - \mu (-1)^{l}) [e_{m + l + 2} e_{2}]. 
  \label{eq:coeff}
\end{align}
For $l$ even, the coefficient is $- l \nu \ne 0$, so we get
$[e_{m + l + 2} e_{2}] = 0$, unless $l \equiv 0 \pmod{p}$. In this case,
we compute
\begin{align*}
  0
  &=
  [e_{m+1} [e_{2} e_{1}^{l-1} e_{2}]]
  \\&=
  [e_{m+1} [e_{2} e_{1}^{l-1}] e_{2}]
  -
  [e_{m+1} e_{2} [e_{2} e_{1}^{l-1}]]
  \\&=
  (\nu - (-1)^{l-1} \nu)
  [e_{m+l+2} e_{2}].
\end{align*}
As $l$ is even here, the coefficient is $2 \nu \ne 0$.

For $l$ odd,  the coefficient in~\eqref{eq:coeff} is $2 \mu  - l \nu$. 
Suppose this vanishes. As $1 \le l <  q - 2$, we have $q > 3$ here, so
that $[e_{m-2} e_{2}] = 0$.  We compute
\begin{align*}
  [e_{m-2} [e_{2} e_{1}^{l+2} e_{2}]]
  &=
  \left(
    \binom{l + 2}{2} \mu - \binom{l + 2}{3} \nu
  \right)
  \cdot
  [e_{m + l + 2} e_{2}]
  \\&=
  \frac{(l + 2) (l + 1)}{6} \mu
  [e_{m + l + 2} e_{2}],
\end{align*}
where  we  have  used the  fact  that  $l  \nu  =  2 \mu$  here.   The
coefficient vanishes when $l + 2  \equiv 0 \pmod{p}$, or $l + 1 \equiv
0 \pmod{p}$. (Except  possibly when $p = 3$, and $l + 1$ or $l + 2$
are divisible by $3$ but not by $9$ -- in this  case the rest of
the discussion  is superfluous. Note  that $l \not\equiv  0 \pmod{p}$
here, otherwise $\mu = \frac{1}{2} l \nu = 0$. Therefore $l \equiv -1,
-2 \pmod{3}$ when $p = 3$, so that $(l+2)(l+1)/6$ is an integer.)

When $l + 2 \equiv 0 \pmod{p}$, we have  $0 = 2 \mu - l \nu = 2 (\mu +
\nu)$, so  that we  are in the  case of~\eqref{eq:end_of_constituent},
with $\mu = \lambda$ and $\nu = -\lambda$ for some $\lambda \ne 0$.
Write $l + 2 = \beta p^{t}$,  with $\beta \not\equiv 0$. It is easy to
see, with an argument we have employed before, that $l + p^{t} < 2 q -
3$, so that $[e_{2} e_{1}^{l+p^{t}} e_{2}] = 0$.
We have then
\begin{align*}
  0
  &=
  [e_{m-p^{t}} [e_{2} e_{1}^{l+p^{t}} e_{2}]]
  \\&=
  [e_{m-p^{t}} [e_{2} e_{1}^{l+p^{t}}] e_{2}]
  \\&=
  \left(
    - \binom{l + p^{t}}{p^{t}} \lambda
    + \binom{l + p^{t}}{p^{t} + 1} (- \lambda)
  \right).
  \cdot
  [e_{m+l+2} e_{2}]
  \\&=
  - \lambda \cdot \binom{l + p^{t} + 1}{p^{t} + 1}
  \cdot
  [e_{m+l+2} e_{2}].
\end{align*}
As
\begin{equation*}
  - \lambda \cdot \binom{l + p^{t} + 1}{p^{t} + 1}
  =
  - \lambda \cdot \binom{\beta p^{t} + p^{t} - 1}{p^{t} + 1}
  \equiv
  \lambda \beta 
  \not\equiv
  0
  \pmod{p},
\end{equation*}
we get $[e_{m+l+2} e_{2}] = 0$.

When $l + 1 \equiv 0 \pmod{p}$, write $l + 1 = \beta p^{t}$,
with $\beta \not\equiv 0 \pmod{p}$. Compute
\begin{align*}
  0
  &=
  [e_{m-p^{t}} [e_{2} e_{1}^{l+p^{t}} e_{2}]]
  \\&=
  \left(
    - \binom{l+p^{t}}{p^{t}} \mu + \binom{l+p^{t}}{p^{t} + 1} \nu
  \right)
  \cdot
  [e_{m + l + 1} e_{2}]. 
\end{align*}
The coefficient here is, up to a sign,
\begin{math}
  \beta (\mu +  \nu).
\end{math}
This cannot vanish, otherwise the two relations $\mu + \nu = 0$ and $0
=  2 \mu  - l  \nu =  2 \mu  + \nu$  would yield  $\mu =  \nu =  0$, a
contradiction.

We now provide the induction step for our assumption that all
constituents are of the form $(\lambda, - \lambda)$.

Suppose first the following constituent is of length $q$. Let
\begin{equation*}
  [e_{m+q} e_{2}] = \mu e_{m+q+2},
  \qquad\text{and}\qquad
  [e_{m+q+1} e_{2}] = \nu e_{m+q+3},
\end{equation*}
We have
\begin{align*}
  0
  &=
  [e_{m-1}, [e_{2 q - 1} e_{2}] - e_{2 q + 1}]
  \\&=
  [e_{m-1}, [e_{2} e_{1}^{2 q - 3} e_{2}]] - [e_{m-1},  e_{2 q + 1}]
  \\&=
  [e_{m-1} [e_{2} e_{1}^{2 q - 3}] e_{2}]
  - [e_{m-1} e_{2} [e_{2} e_{1}^{2 q - 3}]]
  - [e_{m-1} [e_{2} e_{1}^{2 q - 1}]]
\end{align*}
The second term vanishes because $[e_{m-1} e_{2}] = 0$. The first term
is a  multiple of $[e_{m  + 2 q  - 2} e_{2}]  = [e_{(m + q) + q  - 2}
e_{2}]$. If this is non-zero, it exhibits a constituent of length $q
- 2$ or $q - 1$, whereas we have shown $q$ to be a lower bound for the
length of a constituent. Therefore the first term also vanishes.

We are left with
\begin{align*}
  0
  &=
  [e_{m-1} [e_{2} e_{1}^{2 q - 1}]]
  \\&\uptoasign
  (-1)^{1} \binom{2 q - 1}{1} [e_{m} e_{2} e_{1}^{2 q - 2}]
  +
  (-1)^{2} \binom{2 q - 1}{2} [e_{m+1} e_{2} e_{1}^{2 q - 3}]
  \\&\phantom{\uptoasign\ }+
  (-1)^{1 + q} \binom{2 q - 1}{1 + q} 
  [e_{m+q} e_{2} e_{1}^{q - 2}]
  +
  (-1)^{1 + q + 1} \binom{2 q - 1}{1 + q + 1} 
  [e_{m+q+1} e_{2} e_{1}^{q - 3}]
\end{align*}
Now the first two binomial coefficients readily evaluate to $1$, while
for the last two we have, for $l \ge q$,
\begin{align*}
  (-1)^{1 + l} \binom{2 q - 1}{1 + l}
  &=
  (-1)^{l+1} 
  \begin{pmatrix}
    q & + & q - 1\\
    q & + & l - q + 1
  \end{pmatrix}
  \\&\equiv
  - (-1)^{l - q + 1} \binom{q - 1}{l - q  + 1}
  \\&= 
  -1.
\end{align*}
We obtain
\begin{equation*}
  0
  =
  \left(
    \lambda - \lambda - \mu - \nu
  \right)
  \cdot
  e_{m + 2 q},
\end{equation*}
so that $\nu = - \mu$, as requested.

Suppose now the next constituent has length $l > q$, so that in
particular
\begin{equation*}
  [e_{m + 2} e_{2}] = [e_{m + 3} e_{2}] = \dots = [e_{m + q} e_{2}] = 0.
\end{equation*}
We first extend this to show $[e_{m + q + 1} e_{2}] = 0$, so that $l > q +
1$. This follows from
\begin{align*}
  0
  &=
  [e_{m} [e_{2} e_{1}^{q-1} e_{2}]]
  \\&=
  [e_{m} [e_{2} e_{1}^{q-1}] e_{2}]
  -
  [e_{m} e_{2} [e_{2} e_{1}^{q-1}]]
  \\&=
  \left( \lambda - \lambda - \lambda  \right)
  \cdot
  [e_{m + q + 1} e_{2}].
\end{align*}
Suppose now
$[e_{m + l} e_{2}] = \mu e_{m + l + 2}$ and $[e_{m + l + 1} e_{2}] = \nu
e_{m + l + 3}$. We compute
\begin{equation*}
  0
  =
  [e_{m + l - q}, [e_{2 q - 1} e_{2}] - e_{2 q + 1}]
  =
  [e_{m + l - q} [e_{2 q - 1} e_{2}] ]
  -
  [e_{m + l - q} e_{2 q + 1}].
\end{equation*}
Keeping  in mind  that $m  + l  -  q \ge  m +  2$, the  first term  is
immediately seen to vanish.
We are left with
\begin{align*}
  0
  &=
  [e_{m + l - q} e_{2 q + 1}]
  \\&=
  [e_{m + l - q}  [ e_{2} e_{1}^{2 q - 1}]]
  \\&=
  (-1)^{q} \binom{2 q - 1}{q} 
  [e_{m + l - q} e_{2} e_{1}^{q - 1}]
  +
  (-1)^{q + 1} \binom{2 q - 1}{q+ 1}
  [e_{m + l - q + 1} e_{2} e_{1}^{q - 2}]
  \\&=
  ( - \mu - \nu ) \cdot e_{m + l + 1}.
\end{align*}
In  this case,  too,  we obtain  $\nu  = -  \mu$.  This completes  the
induction step.

\section{First constituent of length $q$}

\label{sec:extra}

Let $q$ be a power of $p$ ($q>3$), and suppose that $[e_i,e_2]=0$ for $%
i=3,4,\ldots ,q-1$, and that $[e_q,e_2]\ne 0$. By scaling $e_2$ we may
suppose that $[e_q,e_2]=e_{q+2}$. We show that there is a unique infinite
dimensional Lie algebra $L$ of type 2 satisfying this condition. The Lie
algebra $L$ is defined by the following:

\begin{itemize}
\item  $[e_i,e_2]=0$ for $i=3,4,\ldots ,q-1$,

\item  $[e_q,e_2]=e_{q+2}$, $[e_{q+1},e_2]=-\frac 12e_{q+3}$,

\item  $[e_{kq},e_2]=\frac 12e_{kq+2}$, $[e_{kq+1},e_2]=-\frac 12e_{kq+3}$
for $k=2,3,\ldots $,

\item  $[e_k,e_2]=0$ for $k>q+1$ unless $k=0 \pmod{q}$ or $k=1 \pmod{q}$.
\end{itemize}

Note that in this Lie algebra, if $m>q$ and $n\ge 1$ then 
\[
\lbrack [e_m,e_n,e_1^q]=[e_m,e_1^q,e_n]
\]
so that 
\[
\lbrack [e_m,e_{n+q}]=[e_m,[e_n,e_1^q]]=0
\]
It follows that if $m,n>q$ then $[e_m,e_n]=0$, so that the Lie algebra is
soluble. We give a construction of $L$ in Section 8, and we make use of the
existence of $L$ in the following way. In $L$ we have $[e_n,e_2]=\mu
_ne_{n+2}$ for $n>2$, where $\mu _n=0,1,\frac 12,$ or $-\frac 12$ as
described above. Suppose that we have a Lie algebra $M$ of type 2, where $M$
is spanned by $\{e_i\,|\,i\geq 1\}$, with $[e_i,e_1]=e_{i+1}$ for $i>1$ and $%
[e_n,e_2]=\mu _ne_{n+2}$ for $2<n<2m-2$. Then the relation $[e_m,e_m]=0$ gives
\begin{eqnarray*}
0 &=&[e_m,[e_2,e_1^{m-2}]] \\
&=&\sum_{k=0}^{m-2}(-1)^k\binom{m-2}k[e_m,e_1^k,e_2,e_1^{m-2-k}] \\
&=&\sum_{k=0}^{m-3}(-1)^k\binom{m-2}k\mu
_{m+k}e_{2m}+(-1)^{m-2}[e_{2m-2},e_2].
\end{eqnarray*}
So $[e_{2m-2},e_2]=\mu e_{2m}$ for some $\mu $ which is uniquely determined
by $\{\mu _k\,|\,m\leq k<2m-2\}$. The existence of $L$ implies that $\mu
=\mu _{2m-2}$.

So we assume that $[e_i,e_2]=0$ for $i=3,4,\ldots ,q-1$, and that $%
[e_q,e_2]=e_{q+2}$. The argument just given implies that

\[
\lbrack e_{q+1},e_2]=\mu _{q+1}e_{q+3}=-\frac 12e_{q+3}.
\]
Since $q>3$, $[e_1,e_2,e_2]=0$, and so 
\[
0=[e_{q-1},[e_1,e_2,e_2]]=[e_q,e_2,e_2]=[e_{q+2},e_2].
\]
It follows that $[e_{q+3},e_2]=\mu _{q+3}e_{q+5}=0.$

We now show by induction that $[e_k,e_2]=0$ for $k=q+2,q+3,\ldots ,2q-1$. We
have established the cases $k=q+2$ and $q+3$. So suppose that $q+3<m<2q$,
and suppose that $[e_k,e_2]=0$ for $k=q+2,q+3,\ldots ,m-1$.

Using the argument above, it is only necessary to consider the case when $m$
is odd. Then 
\begin{eqnarray*}
0 &=&[e_{q+1},[e_2,e_1^{m-q-3},e_2]] \\
&=&[e_{q+1},[e_2,e_1^{m-q-3}],e_2]-[e_{q+1},e_2,[e_2,e_1^{m-q-3}]] \\
&=&[e_{q+1},e_2,e_1^{m-q-3},e_2]-(-1)^{m-q-3}[e_{q+1},e_2,e_1^{m-q-3},e_2] \\
&=&-[e_m,e_2].
\end{eqnarray*}
So $[e_k,e_2]=0$ for $k=q+2,q+3,\ldots ,2q-1$, as claimed. Also 
\[
\lbrack e_{2q},e_2]=\mu _{2q}e_{2q+2}=\frac 12e_{2q+2}.
\]

The equations obtained so far leave $[e_{2q+1},e_2]$ undetermined, and so we
suppose that 
\[
\lbrack e_{2q+1},e_2]=\lambda e_{2q+3},
\]
for some $\lambda $. We will show below that $\lambda $ must equal $-\frac 12
$ or $-\frac 14$, but first we show that $[e_k,e_2]=0$ for $2q+1<k<3q$. It
is convenient to subdivide the proof of this into the case when $\lambda =0$
and the case when $\lambda \neq 0$.

First consider the case when $\lambda =0$.

The equation $[[e_2,e_1^q],[e_2,e_1^q]]=0$ gives 
\[
\lbrack e_2,e_1^q,e_2,e_1^q]=[e_2,e_1^{2q},e_2], 
\]
which implies that $[e_{2q+2},e_2]=0$. And 
\[
\lbrack e_{2q},[e_1,e_2,e_2]]=0 
\]
gives 
\[
-[e_{2q+3},e_2]=0, 
\]

So we assume that $2q+3<m<3q$, and that $[e_k,e_2]=0$ for $2q<k<m$. If $m$
is odd then 
\[
0=[e_{2q},[e_2,e_1^{m-2q-2},e_2]]=[e_m,e_2]. 
\]
If $m$ is even and $m<3q-1$ then $[e_{2q-1},[e_2,e_1^{m-2q-1},e_2]]=0$ gives 
\begin{equation}
(m-2q-1)[e_m,e_2]=0.  \label{eq1}
\end{equation}
Also if $2q+3<m<3q$ then $[e_{m-q},e_2]=0$, and so the equation 
\[
\lbrack e_{m-q},[e_2,e_1^{q-2},e_2]]=[e_{m-q},[e_2,e_1^q]] 
\]
gives 
\begin{equation}
((3q-m+1)\frac 12+1)[e_m,e_2]=0.  \label{eq3}
\end{equation}
\From (\ref{eq1}) we see that if $m$ is even and $2q+3<m<3q-1$ then $%
[e_m,e_2]=0$ unless $m=1 \pmod{p}$. But (\ref{eq3}) shows that $[e_m,e_2]=0$
in the case when $m=1 \pmod{p}$, as well as in the case when $m=3q-1$. So $%
[e_k,e_2]=0$ for $2q+1<k<3q$ in the case when $\lambda =0$.

So suppose that $\lambda \neq 0$. As above, we want to show that $%
[e_k,e_2]=0 $ for $2q+1<k<3q$. Since we need the following argument several
times, it is convenient to put it in the form of a lemma.

\begin{lemma}
\label{mrvl1}Let $t\geq 1$ and let $q=p^s>3$. Suppose that $[e_k,e_2]=0$ for 
$1<k<2tq$ unless $k=0 \pmod{q}$ or $k=1 \pmod{q}$, that $[e_{2tq},e_2]=\alpha
e_{2tq+2}$ for some $\alpha \neq 0$, and that $[e_{2tq+1},e_2]=\lambda
e_{2tq+3}$ for some $\lambda \neq 0$. Then $[e_{2tq+k},e_2]=0$ for $1<k<q$.
\end{lemma}

\begin{proof}
The case $k=2$ follows from 
\[
0=[e_{2tq-1},[e_1,e_2,e_2]]=\alpha [e_{2tq+2},e_2]. 
\]

Now suppose by induction that $m$ is odd, that $3\leq m\leq q-2$, and that $%
[e_{2tq+k},e_2]=0$ for all $k$ such that $1<k<m$. We show that $%
[e_{2tq+m},e_2]=[e_{2tq+m+1},e_2]=0$, and this establishes the lemma by
induction on (odd) $m$.

First we have 
\begin{equation}
0=[e_{2tq+1},[e_2,e_1^{m-2},e_2]]=2\lambda [e_{2tq+m+1},e_2]-\lambda
(m-2)[e_{2tq+m},e_2,e_1].  \label{eq21}
\end{equation}
If we let $d=\frac{m+3}2$, we also have 
\begin{eqnarray*}
0 &=&[e_{tq+d},e_{tq+d}] \\
&=&[e_{tq+d},[e_d,e_1^{tq}]] \\
&=&\sum_{r=0}^t(-1)^r\binom tr[e_{(t+r)q+d},e_d,e_1^{(t-r)q}].
\end{eqnarray*}
Now our hypotheses imply that $[e_{(t+r)q+d},e_d]=0$ if $r<t$. So this
equation implies that $[e_{2tq+d},e_d]=0$ also. Since $e_d=[e_2,e_1^{d-2}]$
this gives 
\[
\sum_{r=0}^{d-2}(-1)^r\binom{d-2}r[e_{2tq+d+r},e_2,e_1^{d-2-r}]=0. 
\]
But our inductive hypothesis implies that $[e_{2tq+d+r},e_2]=0$ for $r<d-3$.
So we obtain 
\begin{equation}
\lbrack e_{2tq+m+1},e_2]=(d-2)[e_{2tq+m},e_2,e_1].  \label{eq22}
\end{equation}
Since $d-2=\frac{m-1}2$ (\ref{eq21}) and (\ref{eq22}) imply that $%
[e_{2tq+m+1},e_2]=[e_{2tq+m},e_2,e_1]=0$, which also implies that $%
[e_{2tq+m},e_2]=0$. This completes the proof of the lemma.
\end{proof}

So $[e_k,e_2]=0$ for $2q+1<k<3q$, whatever the value of $\lambda $.

Now consider the equation 
\[
\lbrack e_{2q},[e_2,e_1^{q-2},e_2]]=[e_{2q},[e_2,e_1^q]]. 
\]
This gives 
\begin{eqnarray*}
0
&=&[e_{2q},[e_2,e_1^{q-2}],e_2]-[e_{2q},e_2,[e_2,e_1^{q-2}]]-[e_{2q},e_2,e_1^q]+[e_{2q},e_1^q,e_2]
\\
&=&\frac 12[e_{3q},e_2]+2\lambda [e_{3q},e_2]+\frac 12[e_{3q},e_2]-\frac
12e_{3q+2}+[e_{3q},e_2] \\
&=&(2+2\lambda )[e_{3q},e_2]-\frac 12e_{3q+2},
\end{eqnarray*}
which implies that 
\begin{equation}
e_{3q+3}=(4+4\lambda )[e_{3q},e_2,e_1]  \label{eq5}
\end{equation}
And 
\[
\lbrack e_{2q+1},[e_2,e_1^{q-2},e_2]]=[e_{2q+1},[e_2,e_1^q]] 
\]
gives 
\begin{equation}
2\lambda [e_{3q+1},e_2]-\lambda (q-2)[e_{3q},e_2,e_1]=\lambda
e_{3q+3}-[e_{3q+1},e_2].  \label{eq6}
\end{equation}
In addition, the equation 
\[
\lbrack [e_2,e_1^{\frac{3q-1}2}],[e_2,e_1^{\frac{3q-1}2}]]=0 
\]
gives 
\begin{eqnarray}
0 &=&\left( (-1)^{\frac{q-3}2}\binom{\frac{3q-1}2}{\frac{q-3}2}\frac
12+(-1)^{\frac{q-1}2}\binom{\frac{3q-1}2}{\frac{q-1}2}\lambda \right)
e_{3q+3}  \label{eq7} \\
&&+(-1)^{\frac{3q-3}2}\frac{3q-1}2[e_{3q},e_2,e_1]+(-1)^{\frac{3q-1}%
2}[e_{3q+1},e_2]  \nonumber
\end{eqnarray}
Now 
\[
\binom{\frac{3q-1}2}{\frac{q-3}2}=\binom{q+\frac{q-1}2}{\frac{q-1}2-1}=\frac{%
q-1}2\pmod{p}, 
\]
and 
\[
\binom{\frac{3q-1}2}{\frac{q-1}2}=\binom{q+\frac{q-1}2}{\frac{q-1}2}=1\text{
mod }p, 
\]
So (\ref{eq7}) gives 
\[
(-1)^{\frac{q-1}2}\left( \frac 14+\lambda \right) e_{3q+3}-(-1)^{\frac{3q-3}%
2}\frac 12[e_{3q},e_2,e_1]+(-1)^{\frac{3q-1}2}[e_{3q+1},e_2]=0, 
\]
which implies that 
\begin{equation}
\left( \frac 14+\lambda \right) e_{3q+3}-\frac
12[e_{3q},e_2,e_1]-[e_{3q+1},e_2]=0.  \label{eq8}
\end{equation}
\From (\ref{eq5}) and (\ref{eq8}) we obtain 
\[
\lbrack e_{3q+1},e_2]=((4+4\lambda )(\frac 14+\lambda )-\frac
12)[e_{3q},e_2,e_1]. 
\]
So (\ref{eq6}) gives 
\[
(2\lambda +1)((4+4\lambda )(\frac 14+\lambda )-\frac 12)+2\lambda =\lambda
(4+4\lambda ), 
\]
which implies that $\lambda =-\frac 12$ or $-\frac 14$.

If $\lambda =-\frac 12$ then $[e_{3q},e_2]=\frac 12e_{3q+2}$, and $%
[e_{3q+1},e_2]=-\frac 12e_{3q+3}$. If $\lambda =-\frac 14$ and $p=3$, then (%
\ref{eq5}) gives $e_{3q+3}=0$, so $\lambda =-\frac 12$ is the only
possibility when $p=3$. If $\lambda =-\frac 14$ and $p\neq 3$, then we have $%
[e_{3q},e_2]=\frac 13e_{3q+2}$ and $[e_{3q+1},e_2]=-\frac 16e_{3q+3}$.

Thus we have established that if $q$ is a power of $p$ ($q>3$), and $%
[e_i,e_2]=0$ for $i=3,4,\ldots ,q-1$, and $[e_q,e_2]=e_{q+2}$, then 
\[
\lbrack e_{q+1},e_2]=-\frac 12e_{q+3}, 
\]
\[
\lbrack e_k,e_2]=0\text{ for }q+1<k<2q, 
\]
\[
\lbrack e_{2q},e_2]=\frac 12e_{2q+2}, 
\]
\[
\lbrack e_{2q+1},e_2]=\lambda e_{2q+3}\text{ where }\lambda =-\frac 12\text{
or }\lambda =-\frac 14, 
\]
\[
\lbrack e_k,e_2]=0\text{ for }2q+1<k<3q, 
\]
\[
\lbrack e_{3q},e_2]=\left\{ 
\begin{array}{l}
\frac 12e_{3q+2}\text{ when }\lambda =-\frac 12 \\ 
\frac 13e_{3q+2}\text{ when }\lambda =-\frac 14
\end{array}
\right. , 
\]
\[
\lbrack e_{3q+1},e_2]=\left\{ 
\begin{array}{l}
-\frac 12e_{3q+3}\text{ when }\lambda =-\frac 12 \\ 
-\frac 16e_{3q+3}\text{ when }\lambda =-\frac 14
\end{array}
\right. . 
\]
Furthermore, the case $\lambda =-\frac 14$ can only arise when $p\neq 3$.

\subsection{Generic step for $\lambda =-\frac 12$.}

We assume that $q$ is a power of $p$ ($q>3$) and we assume that

\begin{itemize}
\item  $[e_i,e_2]=0$ for $i=3,4,\ldots ,q-1$,

\item  $[e_q,e_2]=e_{q+2}$, $[e_{q+1},e_2]=-\frac 12e_{q+3}$,

\item  $[e_{kq},e_2]=\frac 12e_{kq+2}$, $[e_{kq+1},e_2]=-\frac 12e_{kq+3}$
for $k=2,3,\ldots ,2n-1$ ($n\geq 2$),

\item $[e_k,e_2]=0$ for $q+1<k<(2n-1)q$  unless $k=0 \pmod{q}$ or $k=1
  \pmod{q}$.
\end{itemize}

We show that

\begin{itemize}
\item  $[e_{(2n-1)q+k},e_2]=0$ for $1<k<q$,

\item  $[e_{2nq},e_2]=\frac 12e_{2nq+2}$, $[e_{2nq+1},e_2]=\lambda e_{2nq+3}$
where $\lambda =-\frac 12$ or $-\frac 14$,

\item  $[e_{2nq+k},e_2]=0$ for $1<k<q$,

\item  if $\lambda =-\frac 12$ then $[e_{(2n+1)q},e_2]=\frac 12e_{(2n+1)q+2}$
and $[e_{(2n+1)q+1},e_2]=-\frac 12e_{(2n+1)q+3}$,

\item  if $\lambda =-\frac 14$ then $[e_{(2n+1)q},e_2]=\frac 13e_{(2n+1)q+2}$
and $[e_{(2n+1)q+1},e_2]=-\frac 16e_{(2n+1)q+3}$.
\end{itemize}

First we show that $[e_{(2n-1)q+k},e_2]=0$ for $1<k<q$. Since $%
[e_{(2n-2)q+k},e_2]=0$, the equation 
\[
\lbrack e_{(2n-2)q+k},[e_2,e_1^{q-2},e_2]]=[e_{(2n-2)q+k},[e_2,e_1^q]] 
\]
gives 
\[
(-1)^{q-k}\left( \binom{q-2}{q-k}+\binom{q-2}{q-k+1}\right) \frac
12[e_{(2n-1)q+k},e_2]=-[e_{(2n-1)q+k},e_2]. 
\]
This implies that $[e_{(2n-1)q+k},e_2]=0$, since 
%AEC2% Slightly reformatted
\begin{multline*}
  1+(-1)^{q-k}\binom{q-2}{q-k}\frac 12+(-1)^{q-k}\binom{q-2}{q-k+1}\frac 12
  =\\=
  1+(q-k+1)\frac 12-(q-k+2)\frac 12\pmod{p}
  =\\=
  \frac 12\pmod{p}.
\end{multline*}

Next consider the equation 
\[
\lbrack e_{(2n-1)q},[e_2,e_1^{q-2},e_2]]=[e_{(2n-1)q},[e_2,e_1^q]] 
\]
This gives 
\begin{eqnarray*}
&&2[[e_{(2n-1)q},e_2,e_1^{q-2},e_2]+2[[e_{(2n-1)q+1},e_2,e_1^{q-3},e_2] \\
&=&[e_{(2n-1)q},e_2,e_1^q]-[e_{(2n-1)q},e_1^q,e_2]
\end{eqnarray*}
which implies that 
\[
\lbrack e_{2nq},e_2]=\frac 12e_{2nq+2}. 
\]

The equations obtained so far leave $[e_{2nq+1},e_2]$ undetermined, and so
we suppose that 
\[
\lbrack e_{2nq+1},e_2]=\lambda e_{2nq+3}, 
\]
for some $\lambda $. We will show below that $\lambda $ must equal $-\frac
12 $ or $-\frac 14$.

First note that the lemma implies that if $\lambda \neq 0$ then $%
[e_{2nq+k},e_2]=0$ for $1<k<q$. We show that $[e_{2nq+k},e_2]=0$ for $1<k<q$
in the case $\lambda =0$ also. So suppose that $\lambda =0$. 
\[
0=[e_{2nq-1},[e_1,e_2,e_2]]=\frac 12[e_{2nq+2},e_2]. 
\]
Also 
\begin{eqnarray*}
0 &=&[e_{2nq},[e_1,e_2,e_2]] \\
&=&[e_{2nq},e_1,e_2,e_2]-2[e_{2nq},e_2,e_1,e_2]]+[e_{2nq},e_2,e_2,e_1] \\
&=&-[e_{2nq+3},e_2],
\end{eqnarray*}
so $[e_{2nq+3},e_2]=0$.

We assume that $3<m<q$, and that $[e_{2nq+k},e_2]=0$ for $1<k<m$. If $m$ is
odd then 
\[
0=[e_{2nq},[e_2,e_1^{m-2},e_2]]=[e_{2nq+m},e_2]. 
\]
If $m$ is even and $m<q-1$ then $[e_{2nq-1},[e_2,e_1^{m-1},e_2]]=0$ gives 
\begin{equation}
(m-1)[e_{2nq+m},e_2]=0.  \label{eq10}
\end{equation}
Also if $3<m<q$ then $[e_{(2n-1)q+m},e_2]=0$, and so the equation 
\[
\lbrack e_{(2n-1)q+m},[e_2,e_1^{q-2},e_2]]=[e_{(2n-1)q+m},[e_2,e_1^q]] 
\]
gives 
\begin{equation}
((q-m+1)\frac 12+1)[e_{2nq+m},e_2]=0.  \label{eq11}
\end{equation}
\From (\ref{eq10}) we see that if $m$ is even and $3<m<q-1$ then $%
[e_{2nq+m},e_2]=0$ unless $m=1 \pmod{p}$. But (\ref{eq11}) shows that $%
[e_{2nq+m},e_2]=0$ in the case when $m=1 \pmod{p}$, as well as in the case
when $m=q-1$. So $[e_{2nq+k},e_2]=0$ for $1<k<q$ in the case when $\lambda
=0 $, as well as in the case $\lambda \neq 0$.

Now consider the equation 
\[
\lbrack e_{2nq},[e_2,e_1^{q-2},e_2]]=[e_{2nq},[e_2,e_1^q]]. 
\]
This gives 
\begin{equation}
e_{(2n+1)q+3}=(4+4\lambda )[e_{(2n+1)q},e_2,e_1]  \label{eq15}
\end{equation}
in exactly the same way as (\ref{eq5}) was obtained from $%
[e_{2q},[e_2,e_1^{q-2},e_2]]=[e_{2q},[e_2,e_1^q]]$. And 
\[
\lbrack e_{2nq+1},[e_2,e_1^{q-2},e_2]]=[e_{2nq+1},[e_2,e_1^q]] 
\]
gives 
\begin{equation}
2\lambda [e_{(2n+1)q+1},e_2]-\lambda (q-2)[e_{(2n+1)q},e_2,e_1]=\lambda
e_{(2n+1)q+3}-[e_{(2n+1)q+1},e_2].  \label{eq16}
\end{equation}
Now consider the equation 
\begin{equation}
\lbrack [e_{nq+2},e_1^{\frac{q-1}2}],[e_{nq+2},e_1^{\frac{q-1}2}]]=0
\label{eq17}
\end{equation}
If we expand $[[e_{nq+2},e_1^{\frac{q-1}2}],[e_{nq+2},e_1^{\frac{q-1}2}]]$
we obtain a sum of the form 
\[
\sum_{r=\frac{q-1}2}^{q-1}\alpha _r[e_{nq+2+r},e_{nq+2},e_1^{q-1-r}]. 
\]
Now 
\begin{eqnarray*}
&&[e_{nq+2+r},e_{nq+2}] \\
&=&[e_{nq+2+r},[e_2,e_1^{nq}]] \\
&=&\sum_{s=0}^n(-1)^s\binom ns[e_{nq+2+r},e_1^{sq},e_2,e_1^{(n-s)q}].
\end{eqnarray*}
If $r<q-2$, then $[e_{nq+2+r},e_1^{sq},e_2,e_1^{(n-s)q}]=0$ for all $s$. If $%
r=q-2$ then 
\[
\lbrack e_{nq+2+r},e_1^{sq},e_2,e_1^{(n-s)q}]=\frac 12e_{(2n+1)q+2} 
\]
for $s<n$, and 
\[
\lbrack e_{nq+2+r},e_1^{sq},e_2,e_1^{(n-s)q}]=[e_{(2n+1)q},e_2] 
\]
for $s=n$. It follows that if $r=q-2$ then 
\[
\sum_{s=0}^n(-1)^s\binom
ns[e_{nq+2+r},e_1^{sq},e_2,e_1^{(n-s)q}]=(-1)^n(\frac
12e_{(2n+1)q+2}-[e_{(2n+1)q},e_2]). 
\]
Similarly, if $r=q-1$ then 
\[
\sum_{s=0}^n(-1)^s\binom
ns[e_{nq+2+r},e_1^{sq},e_2,e_1^{(n-s)q}]=(-1)^n(\lambda
e_{(2n+1)q+3}-[e_{(2n+1)q+1},e_2]). 
\]

So (\ref{eq17}) gives 
\[
\frac 12(\frac 12e_{(2n+1)q+3}-[e_{(2n+1)q},e_2,e_1])+\lambda
e_{(2n+1)q+3}-[e_{(2n+1)q+1},e_2]=0, 
\]

which implies that 
\begin{equation}
\left( \frac 14+\lambda \right) e_{(2n+1)q+3}-\frac
12[e_{(2n+1)q},e_2,e_1]-[e_{(2n+1)q+1},e_2]=0.  \label{eq18}
\end{equation}
Equations (\ref{eq15}), (\ref{eq16}) and (\ref{eq18}) imply that $\lambda
=-\frac 12$ or $-\frac 14$ in exactly the same way as equations (\ref{eq5}),
(\ref{eq6}) and (\ref{eq8}) do. They similarly imply that if $\lambda
=-\frac 12$ then $[e_{(2n+1)q},e_2]=\frac 12e_{(2n+1)q+2}$, and $%
[e_{(2n+1)q+1},e_2]=-\frac 12e_{(2n+1)q+3}$. If $\lambda =-\frac 14$ and $%
p=3 $, then (\ref{eq15}) gives $e_{(2n+1)q+3}=0$, so $\lambda =-\frac 12$ is
the only possibility when $p=3$. If $\lambda =-\frac 14$ and $p\neq 3$, then
we have $[e_{(2n+1)q},e_2]=\frac 13e_{(2n+1)q+2}$ and $[e_{(2n+1)q+1},e_2]=-%
\frac 16e_{(2n+1)q+3}$.

This establishes the generic step for $\lambda =-\frac 12$.

\subsection{Generic step for $\lambda =-\frac 14$.}

We assume that $q$ is a power of $p$ ($p>3$) and we assume that

\begin{itemize}
\item  $[e_i,e_2]=0$ for $i=3,4,\ldots ,q-1$,

\item  $[e_q,e_2]=e_{q+2}$, $[e_{q+1},e_2]=-\frac 12e_{q+3}$,

\item  $[e_{kq},e_2]=\frac 12e_{kq+2}$, $[e_{kq+1},e_2]=-\frac 12e_{kq+3}$
for $k=2,3,\ldots ,2n-1$ ($n\geq 2$),

\item  $[e_{2nq},e_2]=\frac 12e_{2nq+2}$, $[e_{2nq+1},e_2]=-\frac
14e_{2nq+3} $,

\item  There exists $s$ with $1\leq s<p-2$ such that $[e_k,e_2]=0$ for $%
q+1<k<(2n+s)q$ unless $k=0 \pmod{q}$ or $k=1 \pmod{q}$,

\item  $[e_{(2n+k)q},e_2]=\frac 1{k+2}e_{(2n+k)q+2}$ and $%
[e_{(2n+k)q+1},e_2]=-\frac 1{2(k+2)}e_{(2n+k)q+3}$ for $k=1,2,\ldots ,s$.
\end{itemize}

Note that this situation arises from the case $\lambda =-\frac 14$ of the
last section, with $s=1$.

We show that $[e_{(2n+s)q+k},e_2]=0$ for $1<k<q$. In addition we show that
if $s<p-3$ then $[e_{(2n+s+1)q},e_2]=\frac 1{s+3}e_{(2n+s+1)q+2}$, $%
[e_{(2n+s+1)q+1},e_2]=-\frac 1{2(s+3)}e_{(2n+s+1)q+3}$, and we show that if $%
s=p-3$ then $e_{(2n+s+1)q+2}=0$. This contradiction shows that the case $%
\lambda =-\frac 14$ cannot arise in an infinite dimensional Lie algebra of
type 2.

For the moment we suppose that $s<p-3$.

First we show that $[e_{(2n+s)q+k},e_2]=0$ for $1<k<q$. The case $k=2$
follows from 
\[
0=[e_{(2n+s)q-1},[e_1,e_2,e_2]]=\frac 1{s+2}[e_{(2n+s)q+2},e_2]. 
\]
For $k=3$ we have 
\[
\lbrack e_{(2n+s-1)q+3},[e_2,e_1^{q-2},e_2]]=[e_{(2n+s-1)q+3},[e_2,e_1^q]] 
\]
which gives 
\[
((q-2)\frac 1{s+2}-(q-1)\frac
1{2(s+2)})[e_{(2n+s)q+3},e_2]=-[e_{(2n+s)q+3},e_2], 
\]
and so implies that $[e_{(2n+s)q+3},e_2]=0$ unless $2s+1=0 \pmod{p}$. We also
have 
\begin{eqnarray*}
0 &=&[e_{(2n+s)q},[e_1,e_2,e_2]] \\
&=&(-\frac 1{2(s+2)}-\frac 2{s+2})[e_{(2n+s)q+3},e_2],
\end{eqnarray*}
which implies that $[e_{(2n+s)q+3},e_2]=0$ unless $p=5$. Now if $p=5$ then
our assumption that $1\leq s<p-3$ implies that $s=1$ so that $2s+1\neq 0$
mod 5. So $[e_{(2n+s)q+3},e_2]=0$ in every case.

Now suppose that $3<m<q$ and that $[e_{(2n+s)q+k},e_2]=0$ for all $k$ such
that $1<k<m$. If $m$ is even then 
\[
0=[e_{(2n+s)q+1},[e_2,e_1^{m-3},e_2]]=-\frac 1{s+2}[e_{(2n+s)q+m},e_2]. 
\]
So we may assume that $m$ is odd. In this case we have 
\begin{eqnarray*}
0 &=&[e_{(2n+s)q},[e_2,e_1^{m-2},e_2]] \\
&=&(\frac 2{s+2}-(m-2)\frac{-1}{2(s+2)})[e_{(2n+s)q+m},e_2],
\end{eqnarray*}
so $[e_{(2n+s)q+m},e_2]=0$ unless $m=-2 \pmod{p}$. We also have 
\[
\lbrack e_{(2n+s-1)q+m},[e_2,e_1^{q-2},e_2]]=[e_{(2n+s-1)q+m},[e_2,e_1^q]]. 
\]
This implies that 
\[
((q-m+1)\frac 1{s+2}-(q-m+2)\frac 1{2(s+2)}+1)[e_{(2n+1)q+m},e_2]=0, 
\]
which implies that $[e_{(2n+s)q+m},e_2]=0$ unless $m=2(s+2) \pmod{p}$. Since $%
s<p-3$, $[e_{(2n+s)q+m},e_2]=0$ in every case.

Next consider the equation 
\[
\lbrack e_{(2n+s)q},[e_2,e_1^{q-2},e_2]]=[e_{(2n+s)q},[e_2,e_1^q]]. 
\]
This implies that 
\[
(\frac 2{s+2}-\frac 2{2(s+2)}+1)[e_{(2n+s+1)q},e_2]=\frac
1{s+2}e_{(2n+s+1)q+2}, 
\]
and hence that 
\[
\lbrack e_{(2n+s+1)q},e_2]=\frac 1{s+3}e_{(2n+s+1)q+2}. 
\]
And the equation 
\[
\lbrack e_{(2n+s)q+1},[e_2,e_1^{q-2},e_2]]=[e_{(2n+s)q+1},[e_2,e_1^q]] 
\]
implies that 
\begin{eqnarray*}
&&(-\frac 2{2(s+2)}+1)[e_{(2n+s+1)q+1},e_2]+\frac
1{2(s+2)}(q-2)[e_{(2n+s+1)q},e_2,e_1] \\
&=&-\frac 1{2(s+2)}e_{(2n+s+1)q+3}.
\end{eqnarray*}
So 
\[
\lbrack e_{(2n+s+1)q+1},e_2]=-\frac 1{2(s+3)}e_{(2n+s+1)q+3}. 
\]

Finally we consider the case when $s=p-3$. Let $2t=2n+s=2n+p-3$. Then $%
[e_{2tq},e_2]=-e_{2tq+2}$ and $[e_{2tq+1},e_2]=\frac 12e_{2tq+3}$.

The lemma implies that $[e_{2tq+k},e_2]=0$ for $1<k<q$.

Consider the equation 
\[
\lbrack e_{2tq},[e_2,e_1^{q-2},e_2]]=[e_{2tq},[e_2,e_1^q]]. 
\]
This implies that 
\[
(-2-\frac{q-2}2+1)[e_{(2t+1)q},e_2]=-e_{(2t+1)q+2}, 
\]
and hence that 
\[
e_{(2t+1)q+2}=0. 
\]

\section{The case $[e_{3} e_{2}] \ne 0$ and $p =3$}
\label{sec:extra_3}

Let $L$ be an $\mathbb{N}$-graded  Lie algebra of maximal class over a
field $\F$
%AEC% I have used here and in the rest of the file``F'' instead of
%``K'' for consistency. 
of characteristic 3, where $L$ has basis $\{e_i\,|\,i=1,2,\ldots \}$, with $%
[e_i,e_1]=e_{i+1}$ for $i>1$. We consider the case when $[e_3,e_2]\neq 0$.
By rescaling $e_2$ we may assume that $[e_3,e_2]=e_5$, which implies that $%
[e_4,e_2]=e_6$ but leaves $[e_5,e_2]$ undetermined. We show that for every $%
\lambda \in \F$ there is a unique infinite dimensional soluble algebra $%
L(\lambda )$ of type 2 satisfying these relations, together with the
relation $[e_5,e_2]=\lambda e_7$. The algebra $L(\lambda )$ has basis $%
\{e_i\,|\,i=1,2,\ldots \}$, and satisfies the following relations:

\begin{equation}\label{eq:identities_for_Llambda}
  \begin{cases}{}
    [e_i,e_1]=e_{i+1},& \text{for $i>1$,}\\{}
    [e_3,e_2]=e_5,\\{}
        \begin{aligned}[b]
          &[e_{3k+1},e_2]=e_{3k+3}, [e_{3k+2},e_2]=\lambda e_{3k+4},\\
          &\qquad[e_{3k+3},e_2]=(-1-\lambda )e_{3k+5},
        \end{aligned}& \text{for $k\geq 1$,}\\{}
    [e_k,e_3]=(1-\lambda )e_{k+3},& \text{for $k\geq 4$,}\\{}
    [e_k,e_m]=0,& \text{for $k,m\geq 4$.}
  \end{cases}
\end{equation}

% \begin{itemize}
% \item  $[e_i,e_1]=e_{i+1}$ for $i>1$,

% \item  $[e_3,e_2]=e_5$,

% \item  $[e_{3k+1},e_2]=e_{3k+3}$, $[e_{3k+2},e_2]=\lambda e_{3k+4}$, $%
% [e_{3k+3},e_2]=(-1-\lambda )e_{3k+5}$ for $k\geq 1$,

% \item  $[e_k,e_3]=(1-\lambda )e_{k+3}$ for $k\geq 4$,

% \item  $[e_k,e_m]=0$ for $k,m\geq 4$.
% \end{itemize}

We give a construction of $L(\lambda )$ in
Section~\ref{sec:construction_3}, but in fact it is easy 
to show directly that these relations (together with the relations $%
[e_i,e_i]=0$, $[e_i,e_j]+[e_j,e_i]=0$) imply the Jacobi relations 
\[
\lbrack e_i,e_j,e_k]+[e_j,e_k,e_i]+[e_k,e_i,e_j]=0. 
\]
Note that $L(1)\cong m_2$ and that $L(-1)$ is the analogue for $q=3$ of the
algebra constructed in Section 6.

So we suppose that $L$ has basis $\{e_i\,|\,i=1,2,\ldots \}$, with $%
[e_i,e_1]=e_{i+1}$ for $i>1$, $[e_3,e_2]=e_5$, $[e_4,e_2]=e_6$, $%
[e_5,e_2]=\lambda e_7$. We show that if $n\geq 4$ then $[e_n,e_2]=\mu
_ne_{n+2}$, where 
\[
\mu _n=\left\{ 
\begin{array}{l}
1\text{ if }n=1 \pmod{3} \\ 
\lambda \text{ if }n=2 \pmod{3} \\ 
-1-\lambda \text{ if }n=0 \pmod{3}.
\end{array}
\right. 
\]
The fact that $[e_k,e_3]=(1-\lambda )e_{k+3}$ for $k\geq 4$, and that $%
[e_k,e_m]=0$ for $k,m\geq 4$, follows easily from this.

We will make use of the following argument. Suppose that we have shown that $%
[e_n,e_2]=\mu _ne_{n+2}$ for all $n$ with $4\leq n<2m$. Then the relation $%
[e_{m+1},e_{m+1}]=0$ implies that 
\begin{eqnarray*}
0 &=&[e_{m+1},[e_2,e_1^{m-1}]] \\
&=&\sum_{k=0}^{m-1}(-1)^k\binom{m-1}k[e_{m+1},e_1^k,e_2,e_1^{m-1-k}],
\end{eqnarray*}
and so $[e_{2m},e_2]$ is determined by the values of $[e_n,e_2]$ for $%
m+1\leq n<2m$. So $[e_{2m},e_2]=\mu e_{2m+2}$, for some $\mu $ which is
uniquely determined by $\{\mu _n\,|\,m+1\leq n<2m\}$. But $L(\lambda )$ is a
Lie algebra which satisfies $[e_n,e_2]=\mu _ne_{n+2}$ for all $n\geq 4$. So $%
\mu =\mu _{2m}$. In particular, this argument implies that $%
[e_6,e_2]=(-1-\lambda )e_8$.

Now suppose that $[e_n,e_2]=\mu _ne_{n+2}$ for $4\leq n<m$ for some $m\geq 7$%
. We show that this implies that $[e_m,e_2]=\mu _me_{m+2}$. By the argument
above, we only need to consider the case when $m$ is odd. We use the fact
that $[e_2,e_1^3]+[e_1,e_2,e_2]=0$. So 
%AEC2% Reformatted - gave a too long line
\begin{align*}
  0 
  &=
  [e_{m-3},[e_2,e_1^3]]+[e_{m-3},[e_1,e_2,e_2]] 
  \\&=
  [e_{m-3},e_2,e_1^3]-[e_{m-3},e_1^3,e_2]+
  \\&\phantom{=\ }+
  [e_{m-3},e_1,e_2,e_2]+[e_{m-3},e_2,e_1,e_2]+[e_{m-3},e_2,e_2,e_1]
  \\&=
  \mu _{m-3}(1+\mu _{m-1})e_{m+2}-(1-\mu _{m-2}-\mu _{m-3})[e_m,e_2] 
  \\&=
  \mu _{m-3}(1+\mu _{m-1})e_{m+2}-(1+\mu _{m-1})[e_m,e_2].
\end{align*}
Provided $1+\mu _{m-1}\neq 0$, this gives $[e_m,e_2]=\mu _{m-3}e_{m+2}=\mu
_me_{m+2}$, as required. Note that $1+\mu _{m-1}=0$ can only occur when $%
\lambda =0$ and $m=1 \pmod{3}$, or when $\lambda =-1$ and $m=0 \pmod{3}$. So the
uniqueness of $L(\lambda )$ is established except in the cases when $\lambda
=0$ and $\lambda =-1$. We deal with these two cases separately.

\subsection{The case $\lambda =0$.}

Let $L$ be an $\mathbb{N}$-graded Lie algebra spanned by $\{e_i\,|\,i=1,2,%
\ldots \}$, with $[e_i,e_1]=e_{i+1}$ for $i>1$. Let $[e_3,e_2]=e_5$, $%
[e_4,e_2]=e_6$, $[e_5,e_2]=0$. As above, we suppose that for some $n\geq 1$
we have $[e_{3k+1},e_2]=e_{3k+3}$, $[e_{3k+2},e_2]=0$, $%
[e_{3k+3},e_2]=-e_{3k+5}$ for $1\leq k\leq n$, and we suppose that $%
[e_{3n+4},e_2]=\mu e_{3n+6}$ for some $\mu \neq 1$. As above, we may assume
that $n$ is odd. We prove that $L(0)$ is the unique infinite dimensional
algebra over $\F$ of type 2 satisfying $[e_3,e_2]=e_5$, $[e_5,e_2]=0$ by
showing that this implies that $L$ is nilpotent.

First note that 
%AEC2% Reformatted - gave a too long line
\begin{align*}
  0 
  &=
  [e_{3n+2},[e_2,e_1^3]]+[e_{3n+2},[e_1,e_2,e_2]] 
  \\&=
  [e_{3n+2},e_2,e_1^3]-[e_{3n+2},e_1^3,e_2]+
  \\&\phantom{=\ }+
  [e_{3n+2},e_1,e_2,e_2]+[e_{3n+2},e_2,e_1,e_2]+[e_{3n+2},e_2,e_2,e_1]
  \\&=
  [e_{3n+5},e_2]
\end{align*}

We also have $[e_2,e_1^3,e_2]=0$ which implies that $%
[e_{3n+1},[e_2,e_1^3,e_2]]=0$. This gives $(1+\mu )[e_{3n+6},e_2]=e_{3n+8}$.
If $\mu =-1$ then we have $e_{3n+8}=0$ and $L$ is nilpotent (as claimed). So
we assume that $\mu \neq -1$ and that $[e_{3n+6},e_2]=\frac 1{1+\mu
}e_{3n+8} $. Next, 
\[
\lbrack e_{3n+4},[e_2,e_1,e_2]]=[e_{3n+4},[e_2,e_1^3]] 
\]
implies that $[e_{3n+7},e_2]=\frac \mu {1+\mu }e_{3n+9}$. And since $%
[e_4,e_2]=e_6=[e_3,e_1^3]$ we have 
\[
\lbrack e_{3n+3},[e_4,e_2]]=[e_{3n+3},[e_3,e_1^3]], 
\]
which gives 
\[
\mu +\frac 1{1+\mu }=-1-\mu -\frac 1{1+\mu }+\frac \mu {1+\mu }. 
\]
But this implies that $\mu =0$, and so $[e_{3n+4},e_2]=[e_{3n+5},e_2]=0$, $%
[e_{3n+6},e_2]=e_{3n+8}$, and $[e_{3n+7},e_2]=0$.

Now let $m=\frac{n+1}2$. Then 
\begin{eqnarray*}
0 &=&[[e_3,e_1^{3m}],[e_3,e_1^{3m}]] \\
&=&\sum_{r=0}^m(-1)^r\binom mr[e_{3m+3},e_1^{3r},e_3,e_1^{3(m-r)}] \\
&=&\sum_{r=0}^m(-1)^r\binom mr[e_{3m+3+3r},e_3,e_1^{3(m-r)}].
\end{eqnarray*}
Our inductive hypothesis implies that $[e_{3k},e_3]=e_{3k+3}$ for $m+1\leq
k\leq n$. And $[e_{3n+3},e_3]=-e_{3n+6}$, $[e_{3n+6},e_3]=e_{3n+9}$. So this
equation gives $me_{3n+9}=0$. It follows that $m=0 \pmod{3}$, and hence that $%
n=-1 \pmod{3}$.

Since $n=-1 \pmod{3}$, $n>1$, and so $[e_7,e_2]=e_9$. Hence 
\[
\lbrack e_{3n+1},[e_7,e_2]]=[e_{3n+1},e_9]. 
\]
Since $e_9=[e_3,e_1^3,e_1^3]$ we see that 
%AEC2% Slightly reformatted
\begin{eqnarray*}
[e_{3n+1},e_9]
&=&[e_{3n+1},e_3,e_1^3,e_1^3]+[e_{3n+1},e_1^3,e_3,e_1^3]+[e_{3n+1},e_1^3,e_1^3,e_3]
\\
&=&e_{3n+10}-[e_{3n+8},e_2].
\end{eqnarray*}
And since $[e_7,e_2]=[e_4,e_1^3,e_2]$ we have 
%AEC2% Slightly reformatted
\begin{align*}
  [e_{3n+1},[e_7,e_2]]
  &=
  [e_{3n+1},e_4,e_1^3,e_2]-[e_{3n+1},e_1^3,e_4,e_2]
  \\&\phantom{=\ }-[e_{3n+1},e_2,e_4,e_1^3]+[e_{3n+1},e_2,e_1^3,e_4]
  \\&=
  -e_{3n+10}.
\end{align*}
So $[e_{3n+8},e_2]=-e_{3n+10}$.

To summarize, we may assume that $n$ is odd and $n=-1 \pmod{3}$, and that

\begin{itemize}
\item  $[e_{3k},e_2]=-e_{3k+2}$, $[e_{3k+1},e_2]=e_{3k+3}$, $[e_{3k+2},e_2]=0
$ for $2\leq k\leq n$,

\item  $[e_{3n+3},e_2]=-e_{3n+5}$, $[e_{3n+4},e_2]=[e_{3n+5},e_2]=0$,

\item  $[e_{3n+6},e_2]=e_{3n+8}$, $[e_{3n+7},e_2]=0$, $%
[e_{3n+8},e_2]=-e_{3n+10}$.
\end{itemize}

We let $n=2cq-1$, where $q$ is a power of 3 and where $c$ is coprime to 3.
Then we make a further inductive assumption that for some $t$ with $n+2\leq
t\leq n+q $ we have

\begin{itemize}
\item  $[e_{3k},e_2]=e_{3k+2}$, $[e_{3k+1},e_2]=0$, $[e_{3k+2},e_2]=-e_{3k+4}
$ for $n+2\leq k\leq t$.
\end{itemize}

We show that this implies that $[e_{3t+3},e_2]=e_{3t+5}$, $[e_{3t+4},e_2]=0$%
, $[e_{3t+5},e_2]=-e_{3t+7}$. We have to divide the proof that $%
[e_{3t+3},e_2]=e_{3t+5}$ into two cases depending on whether $t$ is odd or
even.

If $t$ is odd let $m=\frac{t-1}2$. Then, since $e_{3m+4}=[e_4,e_1^{3m}]$,
we see that the equation $[e_{3m+4},e_{3m+4}]=0$ gives 
\[
\sum_{r=0}^m(-1)^r\binom mr[e_{3m+4},e_1^{3r},e_4,e_1^{3(m-r)}]=0. 
\]
Now $[e_{3k+1},e_4]=0$ for $m<k\leq n$ and for $n+1<k<t$, $%
[e_{3n+4},e_4]=e_{3n+8}$, $[e_{3t+1},e_4]=-e_{3t+5}+[e_{3t+3},e_2]$. So we
obtain 
\[
(-1)^{n-m}\binom m{n-m}e_{3t+5}-(-1)^me_{3t+5}+(-1)^m[e_{3t+3},e_2]=0, 
\]
which implies that 
\[
\lbrack e_{3t+3},e_2]=(1+\binom m{n-m})e_{3t+5}. 
\]
Now we can write $t=n+2s$ for some $s$ with $2\leq 2s<q$. So 
\[
\binom m{n-m}=\binom m{2m-n}=\binom{cq+s-1}{2s-1}=0\pmod{3},
\]
and $[e_{3t+3},e_2]=e_{3t+5}$.

Now consider the case when $t$ is even. We have $%
e_{3(t-n)}=[e_{3(t-n)-2},e_2]$, and so 
\[
\lbrack e_{3n+5},[e_{3(t-n)-2},e_2]]=[e_{3n+5},e_{3(t-n)}]. 
\]
Now 
%AEC2% Slightly reformatted
\begin{eqnarray*}
[e_{3n+5},e_{3(t-n)}]
&=&[e_{3n+5},[e_3,e_1^{3(t-n-1)}]] \\
&=&\sum_{r=0}^{t-n-1}(-1)^r\binom{t-n-1}%
r[e_{3n+5},e_1^{3r},e_3,e_1^{3(t-n-1-r)}].
\end{eqnarray*}
Since $[e_{3n+5},e_3]=-e_{3n+8}$, $[e_{3k+5},e_3]=e_{3k+8}$ for $n<k<t-1$, $%
[e_{3t+2},e_3]=-e_{3t+5}-[e_{3t+3},e_2]$, and since $t-n-1$ is even, this
implies that 
\[
\lbrack e_{3n+5},e_{3(t-n)}]=-e_{3t+5}-[e_{3t+3},e_2]. 
\]
Also 
%AEC2% Slightly reformatted
\begin{eqnarray*}
[e_{3n+5},[e_{3(t-n)-2},e_2]]
&=&[e_{3n+5},e_{3(t-n)-2},e_2] \\
&=&[e_{3n+5},[e_4,e_1^{3(t-n-2)}],e_2].
\end{eqnarray*}
Since $[e_{3n+5},e_4]=e_{3n+9}$ and $[e_{3k+5},e_4]=0$ for $n<k\leq t-2$,
this implies that 
\[
\lbrack e_{3n+5},[e_{3(t-n)-2},e_2]]=[e_{3t+3},e_2]. 
\]
So the equation 
\[
\lbrack e_{3n+5},[e_{3(t-n)-2},e_2]]=[e_{3n+5},e_{3(t-n)}] 
\]
implies that $[e_{3t+3},e_2]=e_{3t+5}$.

So $[e_{3t+3},e_2]=e_{3t+5}$ whether $t$ is odd or even.

Next note that the equations 
\begin{eqnarray*}
\lbrack e_{3t+1},[e_2,e_1,e_2]] &=&[e_{3t+1},[e_2,e_1^3]], \\
\lbrack e_{3t+2},[e_2,e_1,e_2]] &=&[e_{3t+2},[e_2,e_1^3]]
\end{eqnarray*}
give $[e_{3t+4},e_2]=0$, $[e_{3t+5},e_2]=-e_{3t+7}$.

So, by induction, we may assume that

\begin{itemize}
\item  $[e_{3k},e_2]=e_{3k+2}$, $[e_{3k+1},e_2]=0$, $[e_{3k+2},e_2]=-e_{3k+4}
$ for $n+2\leq k\leq n+q+1$.
\end{itemize}

Finally, let $m=\frac{n+q}2$. We have 
\begin{eqnarray*}
0 &=&[[e_3,e_1^{3m}],[e_3,e_1^{3m}]] \\
&=&\sum_{r=0}^m(-1)^r\binom mr[e_{3m+3},e_1^{3r},e_3,e_1^{3(m-r)}].
\end{eqnarray*}
We have $[e_{3k},e_3]=e_{3k+3}$ for $m+1\leq k\leq n$ and for $n+1<k\leq
n+q+1$, $[e_{3n+3},e_3]=-e_{3n+6}$. Since $\sum_{r=0}^m(-1)^r\binom mr=0$,
we obtain $\binom mqe_{3n+3q+6}=0$. Since $m=cq+\frac{q-1}2$, $\binom
mq=c\neq 0 \pmod{3}$, and so $e_{3n+3q+6}=0$.

Thus the assumption that $[e_{3n+4},e_2]\neq e_{3n+6}$ implies that $L$ is
nilpotent in every case. This completes our analysis of the case when $%
\lambda =0$.

\subsection{The case $\lambda =-1$.}

Let $L$ be an $\mathbb{N}$-graded Lie algebra spanned by $\{e_i\,|\,i=1,2,%
\ldots \}$, with $[e_i,e_1]=e_{i+1}$ for $i>1$. Let $[e_3,e_2]=e_5$, $%
[e_4,e_2]=e_6$, $[e_5,e_2]=-e_7$. Repeating the argument above, we have $%
[e_6,e_2]=\mu _6e_8=0$ (since 6 is even), $[e_7,e_2]=\mu _7e_9=e_9$ (since $%
7\neq 0 \pmod{3}$), and $[e_8,e_2]=\mu _8e_{10}=-e_{10}$ (since 8 is even). And
we may suppose that for some even $n\geq 2$ we have $[e_{3k},e_2]=0$, $%
[e_{3k+1},e_2]=e_{3k+3}$, $[e_{3k+2},e_2]=-e_{3k+4}$ for $2\leq k\leq n$,
and that $[e_{3n+3},e_2]=\mu e_{3n+6}$ for some $\mu \neq 0$. We prove that $%
L(-1)$ is the unique infinite dimensional algebra over $\F$ of type 2
satisfying $[e_3,e_2]=e_5$, $[e_5,e_2]=-e_7$ by showing that this implies
that $L$ is nilpotent.

The relation 
\[
\lbrack e_{3n+1},[e_2,e_1,e_2]]=[e_{3n+1},[e_2,e_1^3]] 
\]
gives $[e_{3n+4},e_2]=(1+\mu )e_{3n+6}$. And the relation 
\[
\lbrack e_{3n+2},[e_2,e_1,e_2]]=[e_{3n+2},[e_2,e_1^3]] 
\]
gives 
\[
-(1+\mu )[e_{3n+5},e_2]=(1-\mu )e_{3n+7}. 
\]
If $\mu =-1$ then this gives $e_{3n+7}=0$, and $L$ is nilpotent. So we may
suppose that $\mu \neq -1$, and that $[e_{3n+5},e_2]=\frac{\mu -1}{\mu +1}%
e_{3n+7}$.

Since $[e_4,e_2]=e_6=[e_3,e_1^3]$ we obtain 
\[
\lbrack e_{3n+2},[e_4,e_2]]=[e_{3n+2},[e_3,e_1^3]].
\]
This gives 
\[
(1-\mu )[e_{3n+6},e_2]+(\frac{\mu -1}{\mu +1}+\mu
+1)e_{3n+8}=[e_{3n+6},e_2]-(\frac{\mu -1}{\mu +1}+\mu +1)e_{3n+8}.
\]
Since $\mu \neq 0$ we have $[e_{3n+6},e_2]=-\frac \mu {\mu +1}e_{3n+8}$. 

Let $n=2m$. Then, since $e_{3m+4}=[e_4,e_1^{3m}]$, the equation $%
[e_{3m+4},e_{3m+4}]=0$ gives 
\[
\sum_{r=1}^m(-1)^m\binom mr[e_{3m+4},e_1^{3r},e_4,e_1^{3(m-r)}]=0. 
\]
Now $[e_{3k+1},e_4]=0$ for $1\leq k<n$, and $[e_{3n+1},e_4]=\mu e_{3n+5}$, $%
[e_{3n+4},e_4]=\frac{\mu (\mu -1)}{\mu +1}$. So we obtain 
\[
(m\mu -\frac{\mu (\mu -1)}{\mu +1})e_{3n+8}=0. 
\]
So either $e_{3n+8}=0$ (and $L$ is nilpotent), or $m\mu (\mu +1)=\mu (\mu
-1) $. But since $\mu \neq 0$, the only solution of $m\mu (\mu +1)=\mu (\mu
-1)$ is $\mu =1$ and $m=0 \pmod{3}$.

So we may assume that $n=2cq$ where $q$ is a power of 3 and where $c$ is
coprime to 3, and we may assume that

\begin{itemize}
\item  $[e_3,e_2]=e_5$, $[e_4,e_2]=e_6$, $[e_5,e_2]=-e_7$,

\item  $[e_{3k},e_2]=0$, $[e_{3k+1},e_2]=e_{3k+3}$, $[e_{3k+2},e_2]=-e_{3k+4}
$ for $2\leq k\leq n$,

\item  $[e_{3n+3},e_2]=e_{3n+5}$, $[e_{3n+4},e_2]=-e_{3n+6}$, $%
[e_{3n+5},e_2]=0$, $[e_{3n+6},e_2]=e_{3n+8}$.
\end{itemize}

We make the further inductive hypothesis that for some $t$ with $n+1\leq
t<n+q$ we have

\begin{itemize}
\item  $[e_{3k+1},e_2]=-e_{3k+3}$, $[e_{3k+2},e_2]=0$, $%
[e_{3k+3},e_2]=e_{3k+5}$ for $n+1\leq k\leq t$.
\end{itemize}

We show that this implies that $[e_{3t+4},e_2]=-e_{3t+6}$, $[e_{3t+5},e_2]=0$%
, $[e_{3t+6},e_2]=e_{3t+8}$.

The equation 
\[
\lbrack e_{3t+1},[e_2,e_1,e_2]]=[e_{3t+1},[e_2,e_1^3]] 
\]
gives $[e_{3t+4},e_2]=-e_{3t+6}$.

We have to divide the proof that $[e_{3t+5},e_2]=0$ into two cases depending
on whether $t$ is odd or even. First suppose that $t$ is odd and let $m=%
\frac{t+1}2$. Then $e_{3m+2}=[e_2,e_1^{3m}]$ and so the equation $%
[e_{3m+2},e_{3m+2}]=0$ gives 
\[
\sum_{r=0}^m(-1)^r\binom mr[e_{3m+2},e_1^{3r},e_2,e_1^{3(m-r)}]=0. 
\]
Now $[e_{3k+2},e_2]=-e_{3k+4}$ for $m\leq k\leq n$, $[e_{3k+2},e_2]=0$ for $%
n<k\leq t$, and so we obtain 
\[
-\sum_{r=0}^{n-m}(-1)^r\binom mre_{3t+7}+(-1)^m[e_{3t+5},e_2]=0. 
\]
We can write $t=n+2s-1=2cq+2s-1$ where $1\leq s<\frac{q+1}2$ so that $m=cq+s$
and $n-m=cq-s$. So 
\[
\sum_{r=0}^{n-m}(-1)^r\binom mr=\sum_{r=0}^{cq-s}(-1)^r\binom{cq+s}r=\pm
\sum_{r=0}^{2s-1}(-1)^r\binom{cq+s}r. 
\]
But $2s-1<q$ and so $\binom{cq+s}r=0 \pmod{3}$ for $s<r\leq 2s-1$. So, working
modulo 3, 
\[
\sum_{r=0}^{n-m}(-1)^r\binom mr=\pm \sum_{r=0}^s(-1)^r\binom{cq+s}r=\pm
\sum_{r=0}^s(-1)^r\binom sr=0, 
\]
and hence $[e_{3t+5},e_2]=0$.

Next suppose that $t$ is even. The equation

\[
\lbrack e_{3n+2},[e_{3(t-n+1)},e_2]]=0 
\]
gives 
\[
\lbrack e_{3n+2},e_{3(t-n+1)},e_2]+[e_{3n+4},e_{3(t-n+1)}]=0. 
\]
Since $e_{3(t-n+1)}=[e_3,e_1^{3(t-n)}]$, this implies that 
\[
\sum_{r=0}^{t-n}(-1)^r\binom{t-n}r\left(
[e_{3n+2},e_1^{3r},e_3,e_1^{3(t-n-r)},e_2]+[e_{3n+4},e_1^{3r},e_3,e_1^{3(t-n-r)}]\right) =0. 
\]
Now $[e_{3n+2},e_3]=e_{3n+5}$, $[e_{3k+2},e_3]=-e_{3k+5}$ for $n<k\leq t$, $%
[e_{3k+4},e_3]=-e_{3k+7}$ for $n\leq k<t$, and $%
[e_{3t+4},e_3]=-e_{3t+7}-[e_{3t+5},e_2]$. Since $t-n$ is even this equation
implies that $[e_{3t+5},e_2]=0$.

So $[e_{3t+5},e_2]=0$ whether $t$ is odd or even.

Finally 
\[
\lbrack e_{3t+3},[e_2,e_1,e_2]]=[e_{3t+3},[e_2,e_1^3]] 
\]
gives $[e_{3t+6},e_2]=e_{3t+8}$. So we may assume by induction that

\begin{itemize}
\item  $[e_{3k+1},e_2]=-e_{3k+3}$, $[e_{3k+2},e_2]=0$, $%
[e_{3k+3},e_2]=e_{3k+5}$ for $n+1\leq k\leq n+q$.
\end{itemize}

To complete our analysis of case 3 we let $t=\frac{n+q-1}2$, and we consider
the equation 
\begin{eqnarray*}
0 &=&[[e_4,e_1^{3t}],[e_4,e_1^{3t}]] \\
&=&\sum_{r=0}^m(-1)^r\binom tr[e_{3t+4},e_1^{3r},e_4,e_1^{3(t-r)}].
\end{eqnarray*}
Since $[e_{3k+4},e_4]=0$ for $t\leq k<n-1$ and for $n\leq k\leq 2t$, and $%
[e_{3n+1},e_4]=e_{3n+5}$ this implies $\binom tqe_{3n+3q+5}=0$. Since $t=cq+%
\frac{q-1}2$, $\binom tq=c\neq 0 \pmod{3}$, and hence $e_{3n+3q+5}=0$.

Thus the assumption that $[e_{3n+3},e_2]\neq 0$ implies that $L$ is
nilpotent in every case. This completes our analysis of the case when $%
\lambda =-1$.

\section{Constructing the algebra\\
with first constituent of length $q$}

\label{sec:Extra_Construction}

In this section we construct the algebra $L$ with first constituent of
length $q$ which is described in Section 6. If $q=3$, this construction
gives the algebra $L(-1)$ of Section 7.

Let $p$ be an odd prime, and let $q$ be a power of $p$. Let $V$ be a vector
space of dimension $q$ over the field $\F(t)$ of rational functions over the
field $\F$ with $p$ elements. We grade $V$ over the cyclic group of order $q$%
, 
\[
V=\Span{v_{0}}\oplus \Span{v_{1}}\oplus \dots \oplus \Span{v_{q-1}}.
\]
Consider the following endomorphisms $D$ and $E$, of $V$. 
\begin{align*}
E& =\begin{cases} v_{i} \mapsto v_{i+1} & \text{if $i \ne q - 1$}\\ v_{q-1}
\mapsto t v_{0}. \end{cases} \\
D& =\begin{cases} v_{0} \mapsto v_{2}\\ v_{q-1} \mapsto - t v_{1}\\ v_{i}
\mapsto 0 & \text{otherwise}. \end{cases}
\end{align*}
Thus $E$ has weight $1$, and $D$ has weight $2$.

We construct the Lie algebra $A$ spanned by $E$ and $D$ in the endomorphism
algebra of $V$.

Consider $[DE^{q-2}]$, which has weight $q\equiv 0\pmod{q}$. For $0\leq j<q$
we have  
\[
v_j[DE^{q-2}]=\sum_{i=0}^{q-2}(-1)^i\binom{q-2}i v_jE^iDE^{q-2-i}.
\]
If $j>0$ then $v_jE^iD=0$ unless $i=q-j-1$ or $q-j$. For $i=q-j$ we have $%
v_jE^i=tv_0$, and thus 
\[
(-1)^i\binom{q-2}i v_jE^iDE^{q-2-i}=(-1)^{q-j}\binom{q-2}{q-j}tv_j,
\]
while for $i=q-1-j$ we have $v_jE^i=v_{q-1}$, and thus 
\begin{align*}
(-1)^i\binom{q-2}i v_jE^iDE^{q-2-i}& =-(-1)^{q-1-j}\binom{q-2}{q-1-j}tv_j \\
& =(-1)^{q-j}\binom{q-2}{q-1-j}tv_j.
\end{align*}

It follows that 
\[
v_j[DE^{q-2}]=v_j(-1)^{q-j}\binom{q-1}{q-j}tv_j=tv_j.
\]
Similarly (for $0\leq i\leq q-2$) we have $v_0E^iD=0$ unless $i=0$, and so $%
v_0[DE^{q-2}]=tv_0$. So $[DE^{q-2}]=t\cdot 1$ is scalar multiplication by $t$%
. It follows that all the $[DE^i]$, for $0\le i\le q-2$, are non-zero, and
thus linearly independent over $\F$, as they have distinct weights $2,\dots
,q$. We claim that $[DE^iD]=0$ for $0\le i<q-2$. To see this, consider the
associative expansion of $[DE^iD]$, which is a linear combination of
monomials of the form $E^\alpha DE^\beta D$, $DE^\beta DE^\alpha $, with $%
\alpha +\beta =i$. Note that if $E^\alpha DE^\beta D$ or $DE^\beta DE^\alpha 
$ is a monomial which occurs in any of these expansions then $\beta <q-3$.
This is trivially true, except in the expansion of $[DE^{q-3}D]$. However in
the expansion of $[DE^{q-3}D]$, $DE^{q-3}D$ appears twice, \emph{but with
opposite signs}. So it is sufficient to show that if $\beta <q-3$ then $%
v_jDE^\beta D=0$ for all $j$. But $v_jD=0$ unless $j=0$ or $q-1$, and 
\[
v_0DE^\beta D=v_{\beta +2}D=0,
\]
\[
v_{q-1}DE^\beta D=-tv_{\beta +1}D=0,
\]
since $0<\beta +1,\beta +2<q-2$.

Therefore 
\[
A = \Span{E, [D E^{i}] : 0 \le i \le q - 2} 
\]
is $q$-dimensional.

Let us now consider the semidirect product $V+\End(V)$, and in it the Lie
algebra $L$ over $\F$ generated by 
\[
e_1=E,\qquad e_2=-\frac 1{2t}\cdot v_1-\frac 12D.
\]
Recursively define $e_{i+1}=[e_ie_1]$, for $i\ge 2$. Note that for $2\le
i\le q$ we have by induction 
\[
e_i=-\frac 1{2t}v_{i-1}-\frac 12[DE^{i-2}].
\]
In particular for $i=q$ we have 
\[
e_q=-\frac 1{2t}\cdot v_{q-1}-\frac 12[DE^{q-2}]=-\frac 1{2t}v_{q-1}-\frac
t2\cdot 1.
\]
Therefore 
\[
e_{q+1}=-\frac 12v_0,\qquad e_{q+2}=-\frac 12v_1,
\]
and we are in $V$ from now on, and further commutation with $e_1$ and $e_2$
is straightforward. In particular, if $0\leq r<q$ and $k\geq 1$, then $%
e_{kq+r+1}=[e_{q+1}e_1^{(k-1)q+r}]=-\frac 12t^{k-1}v_r$, and 
\[
\lbrack e_{kq+r+1},e_2]=\frac 14t^{k-1}v_rD=\left\{ 
\begin{array}{l}
0\text{ unless }r=0\text{ or }q-1 \\ 
-\frac 12e_{kq+3}\text{ if }r=0 \\ 
\frac 12e_{(k+1)q+2}\text{ if }r=q-1.
\end{array}
\right. 
\]

For $2\le i\le q-1$ we have 
\[
\lbrack e_ie_2]=\frac 14\left( \frac 1tv_{i-1}D-\frac
1tv_1[DE^{i-2}]+[DE^{i-2}D]\right) =0,
\]
because of the above, and the easy fact that $v_1[DE^{i-2}]=0$. And 
\[
\lbrack e_qe_2]=\frac 14\left( \frac 1tv_{q-1}D-\frac
1tv_1[DE^{q-2}]+[DE^{q-2}D]\right) =-\frac 12v_1=e_{q+2}.
\]
So $L$ is of maximal class, graded as we want it to be. We have seen that
the first constituent has length $q$, and that $[e_qe_2]=e_{q+2}$, $%
[e_{q+1},e_2]=-\frac 12e_{q+3}$. For $n>q$ we have $[e_n,e_2]=0$
unless $n$ is congruent to $0$
or $1$ modulo $q$, and for $k\geq 2$ we have $[e_{kq},e_2]=\frac
12e_{kq+2}$,
$%
[e_{kq+1},e_2]=-\frac 12e_{kq+3}$.

\section{Constructing the extra algebras for $q = 3$}
\label{sec:construction_3}

In this section we construct the algebras $L(\lambda)$ of
Section~\ref{sec:extra_3}. These are defined over a field $\F$ of
characteristic $3$, for $\lambda \in \F$.

The construction is similar to the one of the previous section. We
rephrase it here in terms of matrices.

Let $t$ be an indeterminate over $K$. Let $v_{1}, v_{2}, v_{3}$ be the
standard basis of the space of row vectors $K(t)^{3}$.
Consider the $3 \times 3$ matrices over $K(t)$
\begin{equation*}
  E =
  \begin{bmatrix}
     & 1 &   \\
     &   & 1 \\
   t &   &   \\  
  \end{bmatrix},
  \qquad
  D =
  \begin{bmatrix}
            &              & 1  \\
   \lambda t  &              &  \\
            & - (1 + \lambda ) t &   \\  
  \end{bmatrix}.
\end{equation*}
where as usual zero entries are omitted. We have
\begin{equation*}
  [D E] = 
  \begin{bmatrix}
   (1 - \lambda) t  &  &   \\
                &  (1 - \lambda) t &  \\
                &                &  (1 - \lambda) t \\  
  \end{bmatrix},
\end{equation*}
a scalar matrix, so that the Lie algebra spanned by $D$ and $E$ has
dimension $3$. Now consider the block $4 \times 4$ matrices
\begin{equation*}
  e_{1} =
  \begin{bmatrix}
    E & 0\\
    0 & 0\\
  \end{bmatrix},
  \qquad
  e_{2} =
  \begin{bmatrix}
    D & 0\\
    v & 0\\
  \end{bmatrix}.
\end{equation*}
Here
\begin{equation*}
  v = \frac{1}{t} \, v_{2} = [0, \frac{1}{t}, 0] \in K(t)^{3}.
\end{equation*}

Consider the Lie algebra $S$ spanned by $e_{1}$ and $e_{2}$.
We compute
\begin{equation*}
  e_{3} = [e_{2} e_{1}] =
    \begin{bmatrix}
    [D E] & 0\\
    v E    & 0\\
  \end{bmatrix},
  \qquad
  e_{4} = [e_{3} e_{1}] =
  \begin{bmatrix}
    0 & 0\\
    v E^{2}   & 0\\
  \end{bmatrix}.
\end{equation*}
Here $v E = 1/t \cdot v_{3}$, and $v E^{2} = [1, 0, 0] = v_{1}$. If we
define $e_{i+1} = [e_{i} e_{1}]$, for $i \ge 2$, we find thus that for
$i \ge 4$ we have
\begin{equation*}
  e_{i} =
  \begin{bmatrix}
    0 & 0\\
    t^{j} v_{k}   & 0\\
  \end{bmatrix},
\end{equation*}
where $1 \le k \le 3$, and $i = 3 (j+1) + k$. It follows that the
algebra $S$ is infinite-dimensional over $K$, with basis $e_{i}$, for
$i \ge 1$.

We have
\begin{equation*}
  [e_{3} e_{2}] =
  \begin{bmatrix}
    0 & 0\\
    v_{2}   & 0\\
  \end{bmatrix}
  = e_{5}
\end{equation*}
and
\begin{equation*}
  [e_{5} e_{2}] =
  \begin{bmatrix}
    0 & 0\\
    \lambda t v_{1}   & 0\\
  \end{bmatrix}
  = \lambda e_{7}.
\end{equation*}
It     is     now     straightforward     to     see     that     all
identities~\eqref{eq:identities_for_Llambda} are  satisfied in $S$, so
that   $S$   is   isomorphic    to   the   algebra   $L(\lambda)$   of
Section~\ref{sec:extra_3}.

\providecommand{\bysame}{\leavevmode\hbox to3em{\hrulefill}\thinspace}

\end{document}